\documentclass[11pt,leqno]{article}
\normalfont
\renewcommand{\baselinestretch}{1.118}

\usepackage{amsfonts}

\newfont{\eulercursive}{eurm10 at 11pt}
\newcommand{\myl}{\mbox{\eulercursive `}}
\newcommand{\mytranspose}[1]{#1^{\mbox{\scriptsize \textsf T}}}

\newcommand{\QED}{\raisebox{0.5mm}{\fbox{\rule{0mm}{1.5mm}\ }}}

\newcounter{myfn}[page]

\newcommand{\MarsAttacksTheorem}{Theorem 2.1}
\newcommand{\FirstMainTheorems}{Theorems 2.1 and 4.1}

\newcommand{\DynkinDiagramFigure}{Figure 2.1}
\newcommand{\BtwoMarsDemo}{Figure 2.2} 
 
\newcommand{\PosetTheorem}{Theorem 5.2} 
\newcommand{\MainResults}{Theorems 2.1, 4.1, and 5.2}

\newcommand{\LemmaList}{Lemmas 3.2, 3.4, 3.5, 3.7, and 3.8}
\newcommand{\ErikssonTheorems}{Theorems 3.1 and 3.3} 
\newcommand{\StrongConvergenceTheorem}{Theorem 3.1}
\newcommand{\StrongConvergenceCorollary}{Lemma 3.2}
\newcommand{\ComparisonTheorem}{Theorem 3.3}
\newcommand{\ComparisonCorollary}{Lemma 3.4}
\newcommand{\ComparisonResults}{Lemmas 3.2 and 3.4}
\newcommand{\NotMarsFriendlyLemma}{Lemma 3.5}
\newcommand{\NotMarsFriendlyCatalog}{Lemma 3.6}
\newcommand{\NotMarsFriendlyFigure}{Figure 3.1}

\newcommand{\EveryNodeFiredLemma}{Lemma 3.7}
\newcommand{\NewLemmaList}{Lemmas 3.2 and 3.7}
\newcommand{\SubgraphLemma}{Lemma 3.8}
\newcommand{\TwoNodeLemma}{Lemma 3.9}
\newcommand{\FourFamiliesSecondLemma}{Lemma 3.10}

\newcommand{\EGCMLemmaList}{Lemmas 3.2, 3.4, 3.5, 3.7, and 3.8}
\newcommand{\SecondMainResult}{Theorem 4.1} 
\newcommand{\ECoxeterGraphFigure}{Figure 4.1} 
\newcommand{\EGCMLinAlgLemma}{Lemma 4.2}
\newcommand{\EGCMVectorLemma}{Lemma 4.3}
\newcommand{\EGCMInadmissibleProp}{Lemma 4.4}
\newcommand{\EGCMAdmissibleProp}{Lemma 4.5}
\newcommand{\EGCMInadmissibleList}{Lemma 4.6}
\newcommand{\EGCMInadmissibleFigure}{Figure 4.2}
\newcommand{\EGCMAdmissibleLemma}{Lemma 4.7}
\newcommand{\EGCMLemmaListForProof}{Lemmas 4.5 and 4.7}

\newcommand{\PosetTheoremConnectedGraph}{Proposition 5.1}

\newcommand{\ErikssonWordProposition}{Proposition 6.1}
\newcommand{\ProofCorollary}{Theorem 6.2}

\newcommand{\IntroNum}{1}
\newcommand{\FirstMainResultNum}{2}
\newcommand{\FirstProofNum}{3}
\newcommand{\SecondMainResultNum}{4}
\newcommand{\RankedPosetsNum}{5}
\newcommand{\AlgebraClassificationNum}{6}
\newcommand{\RemarksNum}{7}

\newcommand{\myA}{\mbox{\sffamily A}}
\newcommand{\myB}{\mbox{\sffamily B}}
\newcommand{\myC}{\mbox{\sffamily C}}
\newcommand{\myD}{\mbox{\sffamily D}}
\newcommand{\myE}{\mbox{\sffamily E}}
\newcommand{\myF}{\mbox{\sffamily F}}
\newcommand{\myG}{\mbox{\sffamily G}}

\newcommand{\selt}{\mathbf{s}} \newcommand{\telt}{\mathbf{t}}
 
 \newcommand{\xelt}{\mathbf{x}}

\newcommand{\ecolor}{\mathbf{edgecolor}}

\newcommand{\myarrow}[1]{\stackrel{#1}{\rightarrow}}

\newcommand{\CircleInteger}[1]{
\setlength{\unitlength}{0.14cm}
\begin{picture}(3,2) 
\put(1,1){\circle{2}}
\put(0.55,0.6){{\tiny #1}}
\end{picture}
}

\newcommand{\CircleIntegerm}{
\setlength{\unitlength}{0.14cm}
\begin{picture}(3,2) 
\put(1,1){\circle{2}}
\put(0.3,0.6){{\tiny $m$}}
\end{picture}
}

\newcommand{\CircleInfty}{
\setlength{\unitlength}{0.14cm}
\begin{picture}(3,2) 
\put(1,1){\circle{2}}
\put(0.1,0.68){{\tiny $\infty$}}
\end{picture}
}

\newcommand{\AnEGraph}{\setlength{\unitlength}{0.75in}
\begin{picture}(6,0.85)
\put(0,0){\begin{picture}(1,0)
            \put(0,0.1){\circle*{0.075}}
            \put(1,0.1){\circle*{0.075}}
            \put(2,0.1){\circle*{0.075}}
            \put(4,0.1){\circle*{0.075}}
            \put(5,0.1){\circle*{0.075}}
            \put(6,0.1){\circle*{0.075}}
            \put(0,0.1){\line(1,0){2}}
            \multiput(2,0.1)(0.4,0){5}{\line(1,0){0.2}}
            \put(4,0.1){\line(1,0){2}}
           \end{picture}}
\end{picture}}

\newcommand{\BnEGraph}{\setlength{\unitlength}{0.75in}
\begin{picture}(6,0.55)
\put(0,0){\begin{picture}(1,0)
            \put(0,0.1){\circle*{0.075}}
            \put(1,0.1){\circle*{0.075}}
            \put(2,0.1){\circle*{0.075}}
            \put(4,0.1){\circle*{0.075}}
            \put(5,0.1){\circle*{0.075}}
            \put(6,0.1){\circle*{0.075}}
            \put(0,0.1){\line(1,0){2}}
            \multiput(2,0.1)(0.4,0){5}{\line(1,0){0.2}}
            \put(4,0.1){\line(1,0){2}}
            \put(5.3,0.15){\CircleInteger{4}}
           \end{picture}}
\end{picture}}

\newcommand{\DnEGraph}{\setlength{\unitlength}{0.75in}
\begin{picture}(6,0.75)
\put(0,-0.25){\begin{picture}(1,0)
            \put(0,0.35){\circle*{0.075}}
            \put(1,0.35){\circle*{0.075}}
            \put(2,0.35){\circle*{0.075}}
            \put(4,0.35){\circle*{0.075}}
            \put(5,0.35){\circle*{0.075}}
            \put(6,0.1){\circle*{0.075}}
            \put(6,0.6){\circle*{0.075}}
            \put(0,0.35){\line(1,0){2}}
            \multiput(2,0.35)(0.4,0){5}{\line(1,0){0.2}}
            \put(4,0.35){\line(1,0){1}}
            \put(5,0.35){\line(4,1){1}}
            \put(5,0.35){\line(4,-1){1}}
           \end{picture}}
\end{picture}}

\newcommand{\EEightEGraph}{\setlength{\unitlength}{0.75in}
\begin{picture}(6,0.75)
\put(0,-0.25){\begin{picture}(1,0)
            \put(0,0.1){\circle*{0.075}}
            \put(1,0.1){\circle*{0.075}}
            \put(2,0.1){\circle*{0.075}}
            \put(2,0.6){\circle*{0.075}}
            \put(3,0.1){\circle*{0.075}}
            \put(4,0.1){\circle*{0.075}}
            \put(5,0.1){\circle*{0.075}}
            \put(6,0.1){\circle*{0.075}}
            \put(0,0.1){\line(1,0){6}}
            \put(2,0.1){\line(0,1){0.5}}
           \end{picture}}
\end{picture}}

\newcommand{\ESevenEGraph}{\setlength{\unitlength}{0.75in}
\begin{picture}(6,0.75)
\put(0,-0.25){\begin{picture}(1,0)
            \put(0,0.1){\circle*{0.075}}
            \put(1,0.1){\circle*{0.075}}
            \put(2,0.1){\circle*{0.075}}
            \put(2,0.6){\circle*{0.075}}
            \put(3,0.1){\circle*{0.075}}
            \put(4,0.1){\circle*{0.075}}
            \put(5,0.1){\circle*{0.075}}
            \put(0,0.1){\line(1,0){5}}
            \put(2,0.1){\line(0,1){0.5}}
           \end{picture}}
\end{picture}}

\newcommand{\ESixEGraph}{\setlength{\unitlength}{0.75in}
\begin{picture}(6,0.75)
\put(0,-0.25){\begin{picture}(1,0)
            \put(0,0.1){\circle*{0.075}}
            \put(1,0.1){\circle*{0.075}}
            \put(2,0.1){\circle*{0.075}}
            \put(2,0.6){\circle*{0.075}}
            \put(3,0.1){\circle*{0.075}}
            \put(4,0.1){\circle*{0.075}}
            \put(0,0.1){\line(1,0){4}}
            \put(2,0.1){\line(0,1){0.5}}
           \end{picture}}
\end{picture}}

\newcommand{\FFourEGraph}{\setlength{\unitlength}{0.75in}
\begin{picture}(1,0.85)
\put(0,0){\begin{picture}(1,0)
            \put(0,0.1){\circle*{0.075}}
            \put(1,0.1){\circle*{0.075}}
            \put(2,0.1){\circle*{0.075}}
            \put(3,0.1){\circle*{0.075}}
            \put(0,0.1){\line(1,0){3}}
            \put(1.3,0.15){\CircleInteger{4}}
            \end{picture}}
\end{picture}}

\newcommand{\HFourEGraph}{\setlength{\unitlength}{0.75in}
\begin{picture}(1,0.55)
\put(0,0){\begin{picture}(1,0)
            \put(0,0.1){\circle*{0.075}}
            \put(1,0.1){\circle*{0.075}}
            \put(2,0.1){\circle*{0.075}}
            \put(3,0.1){\circle*{0.075}}
            \put(0,0.1){\line(1,0){3}}
            \put(0.3,0.15){\CircleInteger{5}}
            \end{picture}}
\end{picture}}

\newcommand{\HThreeEGraph}{\setlength{\unitlength}{0.75in}
\begin{picture}(1,0.55)
\put(0,0){\begin{picture}(1,0)
            \put(0,0.1){\circle*{0.075}}
            \put(1,0.1){\circle*{0.075}}
            \put(2,0.1){\circle*{0.075}}
            \put(0,0.1){\line(1,0){2}}
            \put(0.3,0.15){\CircleInteger{5}}
            \end{picture}}
\end{picture}}

\newcommand{\ITwoEGraph}{\setlength{\unitlength}{0.75in}
\begin{picture}(1,0.55)
\put(0,0){\begin{picture}(1,0)
            \put(0,0.1){\circle*{0.075}}
            \put(1,0.1){\circle*{0.075}}
            \put(0,0.1){\line(1,0){1}}
            \put(0.3,0.15){\CircleIntegerm}
            \end{picture}}
\end{picture}}

\newcommand{\BTwoGraphForFigure}{\setlength{\unitlength}{0.6in}
\begin{picture}(1,0.4)
\put(0,0){\begin{picture}(1,0)
            \put(0,0.1){\circle*{0.075}}
            \put(-0.11,-0.05){\scriptsize $\gamma_{1}$}
            \put(1,0.1){\circle*{0.075}}
            \put(1,-0.05){\scriptsize $\gamma_{2}$}
            \put(0,0.1){\line(1,0){1}}
            \put(0.2,0.1){\vector(1,0){0.1}}
            \put(0.8,0.1){\vector(-1,0){0.1}}
            \put(0.7,0.1){\vector(-1,0){0.1}}
            \end{picture}}
\end{picture}}

\newcommand{\GTwoMarsOK}{\setlength{\unitlength}{0.75in}
\begin{picture}(1,0.55)
\put(0,0){\begin{picture}(1,0)
            \put(0,0.1){\circle*{0.075}}
            \put(1,0.1){\circle*{0.075}}
            \put(0,0.1){\line(1,0){1}}
            \put(0.2,0.1){\vector(1,0){0.1}}
            \put(0.8,0.1){\vector(-1,0){0.1}}
            \put(0.7,0.1){\vector(-1,0){0.1}}
            \put(0.6,0.1){\vector(-1,0){0.1}}
            \end{picture}}
\end{picture}}

\newcommand{\FFourMarsOK}{\setlength{\unitlength}{0.75in}
\begin{picture}(1,0.85)
\put(0,0){\begin{picture}(1,0)
            \put(0,0.1){\circle*{0.075}}
            \put(1,0.1){\circle*{0.075}}
            \put(2,0.1){\circle*{0.075}}
            \put(3,0.1){\circle*{0.075}}
            \put(0,0.1){\line(1,0){3}}
            \put(1.2,0.1){\vector(1,0){0.1}}
            \put(1.3,0.1){\vector(1,0){0.1}}
            \put(1.8,0.1){\vector(-1,0){0.1}}
            \end{picture}}
\end{picture}}

\newcommand{\CnMarsOK}{\setlength{\unitlength}{0.75in}
\begin{picture}(6,0.55)
\put(0,0){\begin{picture}(1,0)
            \put(0,0.1){\circle*{0.075}}
            \put(1,0.1){\circle*{0.075}}
            \put(2,0.1){\circle*{0.075}}
            \put(4,0.1){\circle*{0.075}}
            \put(5,0.1){\circle*{0.075}}
            \put(6,0.1){\circle*{0.075}}
            \put(0,0.1){\line(1,0){2}}
            \multiput(2,0.1)(0.4,0){5}{\line(1,0){0.2}}
            \put(4,0.1){\line(1,0){2}}
            \put(5.8,0.1){\vector(-1,0){0.1}}
            \put(5.7,0.1){\vector(-1,0){0.1}}
            \put(5.2,0.1){\vector(1,0){0.1}}
           \end{picture}}
\end{picture}}

\newcommand{\BnMarsOK}{\setlength{\unitlength}{0.75in}
\begin{picture}(6,0.55)
\put(0,0){\begin{picture}(1,0)
            \put(0,0.1){\circle*{0.075}}
            \put(1,0.1){\circle*{0.075}}
            \put(2,0.1){\circle*{0.075}}
            \put(4,0.1){\circle*{0.075}}
            \put(5,0.1){\circle*{0.075}}
            \put(6,0.1){\circle*{0.075}}
            \put(0,0.1){\line(1,0){2}}
            \multiput(2,0.1)(0.4,0){5}{\line(1,0){0.2}}
            \put(4,0.1){\line(1,0){2}}
            \put(5.8,0.1){\vector(-1,0){0.1}}
            \put(5.2,0.1){\vector(1,0){0.1}}
            \put(5.3,0.1){\vector(1,0){0.1}}
           \end{picture}}
\end{picture}}

\newcommand{\DnMarsOK}{\setlength{\unitlength}{0.75in}
\begin{picture}(6,0.75)
\put(0,-0.25){\begin{picture}(1,0)
            \put(0,0.35){\circle*{0.075}}
            \put(1,0.35){\circle*{0.075}}
            \put(2,0.35){\circle*{0.075}}
            \put(4,0.35){\circle*{0.075}}
            \put(5,0.35){\circle*{0.075}}
            \put(6,0.1){\circle*{0.075}}
            \put(6,0.6){\circle*{0.075}}
            \put(0,0.35){\line(1,0){2}}
            \multiput(2,0.35)(0.4,0){5}{\line(1,0){0.2}}
            \put(4,0.35){\line(1,0){1}}
            \put(5,0.35){\line(4,1){1}}
            \put(5,0.35){\line(4,-1){1}}
           \end{picture}}
\end{picture}}

\newcommand{\AnMarsOK}{\setlength{\unitlength}{0.75in}
\begin{picture}(6,0.85)
\put(0,0){\begin{picture}(1,0)
            \put(0,0.1){\circle*{0.075}}
            \put(1,0.1){\circle*{0.075}}
            \put(2,0.1){\circle*{0.075}}
            \put(4,0.1){\circle*{0.075}}
            \put(5,0.1){\circle*{0.075}}
            \put(6,0.1){\circle*{0.075}}
            \put(0,0.1){\line(1,0){2}}
            \multiput(2,0.1)(0.4,0){5}{\line(1,0){0.2}}
            \put(4,0.1){\line(1,0){2}}
           \end{picture}}
\end{picture}}

\newcommand{\EEightMarsOK}{\setlength{\unitlength}{0.75in}
\begin{picture}(6,0.75)
\put(0,-0.25){\begin{picture}(1,0)
            \put(0,0.1){\circle*{0.075}}
            \put(1,0.1){\circle*{0.075}}
            \put(2,0.1){\circle*{0.075}}
            \put(2,0.6){\circle*{0.075}}
            \put(3,0.1){\circle*{0.075}}
            \put(4,0.1){\circle*{0.075}}
            \put(5,0.1){\circle*{0.075}}
            \put(6,0.1){\circle*{0.075}}
            \put(0,0.1){\line(1,0){6}}
            \put(2,0.1){\line(0,1){0.5}}
           \end{picture}}
\end{picture}}

\newcommand{\ESevenMarsOK}{\setlength{\unitlength}{0.75in}
\begin{picture}(6,0.75)
\put(0,-0.25){\begin{picture}(1,0)
            \put(0,0.1){\circle*{0.075}}
            \put(1,0.1){\circle*{0.075}}
            \put(2,0.1){\circle*{0.075}}
            \put(2,0.6){\circle*{0.075}}
            \put(3,0.1){\circle*{0.075}}
            \put(4,0.1){\circle*{0.075}}
            \put(5,0.1){\circle*{0.075}}
            \put(0,0.1){\line(1,0){5}}
            \put(2,0.1){\line(0,1){0.5}}
           \end{picture}}
\end{picture}}

\newcommand{\ESixMarsOK}{\setlength{\unitlength}{0.75in}
\begin{picture}(6,0.75)
\put(0,-0.25){\begin{picture}(1,0)
            \put(0,0.1){\circle*{0.075}}
            \put(1,0.1){\circle*{0.075}}
            \put(2,0.1){\circle*{0.075}}
            \put(2,0.6){\circle*{0.075}}
            \put(3,0.1){\circle*{0.075}}
            \put(4,0.1){\circle*{0.075}}
            \put(0,0.1){\line(1,0){4}}
            \put(2,0.1){\line(0,1){0.5}}
           \end{picture}}
\end{picture}}

\newcommand{\TwoCitiesGraphWithLabels}{
\setlength{\unitlength}{0.75in}
\begin{picture}(1.65,0.25)
\put(0.25,0){\begin{picture}(1,0)
            \put(0,0.1){\circle*{0.05}}
            \put(-0.20,-0.05){\large $\gamma_{1}$}
            \put(1,0.1){\circle*{0.05}}
            \put(1.05,-0.05){\large $\gamma_{2}$}
            \put(0,0.1){\line(1,0){1}}
            \put(0.2,0.1){\vector(1,0){0.1}}
            \put(0.8,0.1){\vector(-1,0){0.1}}
            \put(0.225,-0.05){\footnotesize $p$}
            \put(0.71,-0.05){\footnotesize $q$}
            \end{picture}}
\end{picture}}

\newcommand{\EGCMGraphCirclem}{\setlength{\unitlength}{0.75in}
\begin{picture}(1.65,0.3)
\put(0.25,0){\begin{picture}(1,0)
            \put(0,0.1){\circle*{0.05}}
            \put(-0.20,-0.05){\large $\gamma_{1}$}
            \put(1,0.1){\circle*{0.05}}
            \put(1.05,-0.05){\large $\gamma_{2}$}
            \put(0,0.1){\line(1,0){1}}
            \put(0.3,0.15){\CircleIntegerm}
            \end{picture}}
\end{picture}}

\newcommand{\ATwoGraphWithLabels}{\setlength{\unitlength}{0.75in}
\begin{picture}(1,0.3)
\put(0.4,0.225){\large $\mbox{\sffamily A}_{2}$}
\put(0,0){\begin{picture}(1,0)
            \put(0,0.1){\circle*{0.05}}
            \put(-0.15,-0.05){\large $\gamma_{1}$}
            \put(1,0.1){\circle*{0.05}}
            \put(1.05,-0.05){\large $\gamma_{2}$}
            \put(0,0.1){\line(1,0){1}}
            \put(0.2,0.1){\vector(1,0){0.1}}
            \put(0.8,0.1){\vector(-1,0){0.1}}
            \end{picture}}
\end{picture}}

\newcommand{\ATwoGraphNoEdgeLabels}{\setlength{\unitlength}{0.75in}
\begin{picture}(1.65,0.25)
\put(0.25,0){\begin{picture}(1,0)
            \put(0,0.1){\circle*{0.05}}
            \put(-0.20,-0.05){\large $\gamma_{1}$}
            \put(1,0.1){\circle*{0.05}}
            \put(1.05,-0.05){\large $\gamma_{2}$}
            \put(0,0.1){\line(1,0){1}}
            \end{picture}}
\end{picture}}

\newcommand{\BTwoGraphWithLabels}{\setlength{\unitlength}{0.75in}
\begin{picture}(1,0.35)
\put(0.4,0.225){\large $\mbox{\sffamily B}_{2}$}
\put(0,0){\begin{picture}(1,0)
            \put(0,0.1){\circle*{0.05}}
            \put(-0.15,-0.05){\large $\gamma_{1}$}
            \put(1,0.1){\circle*{0.05}}
            \put(1.05,-0.05){\large $\gamma_{2}$}
            \put(0,0.1){\line(1,0){1}}
            \put(0.2,0.1){\vector(1,0){0.1}}
            \put(0.8,0.1){\vector(-1,0){0.1}}
            \put(0.7,0.1){\vector(-1,0){0.1}}
            \end{picture}}
\end{picture}}

\newcommand{\GTwoGraphWithLabels}{\setlength{\unitlength}{0.75in}
\begin{picture}(1,0.35)
\put(0.4,0.225){\large $\mbox{\sffamily G}_{2}$}
\put(0,0){\begin{picture}(1,0)
            \put(0,0.1){\circle*{0.05}}
            \put(-0.15,-0.05){\large $\gamma_{1}$}
            \put(1,0.1){\circle*{0.05}}
            \put(1.05,-0.05){\large $\gamma_{2}$}
            \put(0,0.1){\line(1,0){1}}
            \put(0.2,0.1){\vector(1,0){0.1}}
            \put(0.8,0.1){\vector(-1,0){0.1}}
            \put(0.7,0.1){\vector(-1,0){0.1}}
            \put(0.6,0.1){\vector(-1,0){0.1}}
            \end{picture}}
\end{picture}}

\newcommand{\SmallCycles}{\setlength{\unitlength}{0.75in}
\begin{picture}(6.5,3)
\put(-0.5,1.5){\begin{picture}(1,0)
            \put(0,0.6){\circle*{0.075}}
            \put(0.5,0.1){\circle*{0.075}}
            \put(0.5,1.1){\circle*{0.075}}
            \put(0,0.6){\line(1,1){0.5}}
            \put(0,0.6){\line(1,-1){0.5}}
            \put(0.5,0.1){\line(0,1){1}}
            \put(0.5,0.3){\vector(0,1){0.1}}
            \put(0.5,0.9){\vector(0,-1){0.1}}
            \put(0,0.6){\vector(1,1){0.2}}
            \put(0.5,1.1){\vector(-1,-1){0.2}}
            \put(0,0.6){\vector(1,-1){0.2}}
            \put(0.5,0.1){\vector(-1,1){0.2}}
            \put(-0.05,0.8){\footnotesize $q_{1}$}
            \put(0.175,1){\footnotesize $p_{1}$}
            \put(-0.05,0.35){\footnotesize $q_{2}$}
            \put(0.175,0.15){\footnotesize $p_{2}$} 
            \put(0.7,0.6){\footnotesize $p_{1}q_{1} \geq 1, 
            p_{2}q_{2} \geq 1$}
           \end{picture}}
\put(2.3,1.5){\begin{picture}(1,0)
            \put(0,0.6){\circle*{0.075}}
            \put(0.5,0.1){\circle*{0.075}}
            \put(0.5,1.1){\circle*{0.075}}
            \put(0,0.6){\line(1,1){0.5}}
            \put(0,0.6){\line(1,-1){0.5}}
            \put(0.5,0.1){\line(0,1){1}}
            \put(0.5,0.3){\vector(0,1){0.1}}
            \put(0.5,0.9){\vector(0,-1){0.1}}
            \put(0.5,0.8){\vector(0,-1){0.1}}
            \put(0,0.6){\vector(1,1){0.2}}
            \put(0.5,1.1){\vector(-1,-1){0.2}}
            \put(0,0.6){\vector(1,-1){0.2}}
            \put(0.5,0.1){\vector(-1,1){0.2}}
            \put(-0.05,0.8){\footnotesize $q_{1}$}
            \put(0.175,1){\footnotesize $p_{1}$}
            \put(-0.05,0.35){\footnotesize $q_{2}$}
            \put(0.175,0.15){\footnotesize $p_{2}$} 
            \put(0.7,0.6){\footnotesize $p_{1}q_{1} \geq 2, 
            p_{2}q_{2} \geq 2$}
           \end{picture}}
\put(5.15,1.5){\begin{picture}(1,0)
            \put(0,0.6){\circle*{0.075}}
            \put(0.5,0.1){\circle*{0.075}}
            \put(0.5,1.1){\circle*{0.075}}
            \put(0,0.6){\line(1,1){0.5}}
            \put(0,0.6){\line(1,-1){0.5}}
            \put(0.5,0.1){\line(0,1){1}}
            \put(0.5,0.3){\vector(0,1){0.1}}
            \put(0.5,0.9){\vector(0,-1){0.1}}
            \put(0.5,0.8){\vector(0,-1){0.1}}
            \put(0.5,0.7){\vector(0,-1){0.1}}
            \put(0,0.6){\vector(1,1){0.2}}
            \put(0.5,1.1){\vector(-1,-1){0.2}}
            \put(0,0.6){\vector(1,-1){0.2}}
            \put(0.5,0.1){\vector(-1,1){0.2}}
            \put(-0.05,0.8){\footnotesize $q_{1}$}
            \put(0.175,1){\footnotesize $p_{1}$}
            \put(-0.05,0.35){\footnotesize $q_{2}$}
            \put(0.175,0.15){\footnotesize $p_{2}$} 
            \put(0.7,0.6){\footnotesize $p_{1}q_{1} \geq 3$, 
            $p_{2}q_{2} \geq 3$}
           \end{picture}}
\put(0.5,0){\begin{picture}(1,0)
            \put(0,0.6){\circle*{0.075}}
            \put(0.5,0.1){\circle*{0.075}}
            \put(0.5,1.1){\circle*{0.075}}
            \put(1,0.6){\circle*{0.075}}
            \put(0,0.6){\line(1,1){0.5}}
            \put(0,0.6){\line(1,-1){0.5}}
            \put(1,0.6){\line(-1,1){0.5}}
            \put(1,0.6){\line(-1,-1){0.5}}
            \put(0.95,0.65){\vector(-1,1){0.1}}
            \put(0.55,1.05){\vector(1,-1){0.1}}
            \put(0.65,0.95){\vector(1,-1){0.1}}
           \end{picture}}
\put(2,0){\begin{picture}(1,0)
            \put(0,0.6){\circle*{0.075}}
            \put(0.5,0.1){\circle*{0.075}}
            \put(0.5,1.1){\circle*{0.075}}
            \put(1,0.6){\circle*{0.075}}
            \put(0,0.6){\line(1,1){0.5}}
            \put(0,0.6){\line(1,-1){0.5}}
            \put(1,0.6){\line(-1,1){0.5}}
            \put(1,0.6){\line(-1,-1){0.5}}
            \put(0.95,0.65){\vector(-1,1){0.1}}
            \put(0.55,1.05){\vector(1,-1){0.1}}
            \put(0.65,0.95){\vector(1,-1){0.1}}
            \put(0.45,0.15){\vector(-1,1){0.1}}
            \put(0.05,0.55){\vector(1,-1){0.1}}
            \put(0.15,0.45){\vector(1,-1){0.1}}
           \end{picture}}
\put(3.5,0){\begin{picture}(1,0)
            \put(0,0.6){\circle*{0.075}}
            \put(0.5,0.1){\circle*{0.075}}
            \put(0.5,1.1){\circle*{0.075}}
            \put(1,0.6){\circle*{0.075}}
            \put(0,0.6){\line(1,1){0.5}}
            \put(0,0.6){\line(1,-1){0.5}}
            \put(1,0.6){\line(-1,1){0.5}}
            \put(1,0.6){\line(-1,-1){0.5}}
            \put(0.95,0.65){\vector(-1,1){0.1}}
            \put(0.55,1.05){\vector(1,-1){0.1}}
            \put(0.65,0.95){\vector(1,-1){0.1}}
            \put(0.45,0.15){\vector(-1,1){0.1}}
            \put(0.05,0.55){\vector(1,-1){0.1}}
            \put(0.35,0.25){\vector(-1,1){0.1}}
           \end{picture}}
\put(5,0){\begin{picture}(1,0)
            \put(0,0.6){\circle*{0.075}}
            \put(0,0.1){\circle*{0.075}}
            \put(0.5,0.1){\circle*{0.075}}
            \put(0.5,1.1){\circle*{0.075}}
            \put(1,0.6){\circle*{0.075}}
            \put(0,0.6){\line(1,1){0.5}}
            \put(0,0.6){\line(0,-1){0.5}}
            \put(0,0.1){\line(1,0){0.5}}
            \put(1,0.6){\line(-1,1){0.5}}
            \put(1,0.6){\line(-1,-1){0.5}}
            \put(0.95,0.65){\vector(-1,1){0.1}}
            \put(0.55,1.05){\vector(1,-1){0.1}}
            \put(0.65,0.95){\vector(1,-1){0.1}}
           \end{picture}}
\end{picture}}

\begin{document}

\newpage
\setcounter{page}{1} 
\renewcommand{\baselinestretch}{1}

\vspace*{-0.7in}
\hfill {\footnotesize October 29, 2008}

\begin{center}
{\large \bf The numbers game and 
Dynkin diagram classification results}

Robert G.\ Donnelly$^{1}$ and Kimmo Eriksson$^{2}$

\vspace*{-0.05in} 
$^{1}$Department of Mathematics and Statistics, Murray State
University, Murray, KY 42071

\vspace*{-0.05in} 
$^{2}$Department of Mathematics and Physics, M\"{a}lardalen 
University, V\"{a}ster{\aa}s, Sweden
\end{center}

\begin{abstract}
The numbers game is a one-player game played on  
a finite simple 
graph with certain ``amplitudes'' assigned to its edges 
and with an initial assignment of real 
numbers to its nodes.  The moves of the game successively transform 
the numbers at the nodes using the amplitudes in a certain way.  
Combinatorial reasoning is used to
show that those connected 
graphs with negative integer amplitudes for which the numbers game  
meets a certain finiteness requirement are 
precisely the Dynkin diagrams associated with the finite-dimensional 
complex simple Lie algebras. 
This strengthens a result originally due to the second author. 
A more general result is obtained when certain real number 
amplitudes are allowed. The resulting graphs are in families, each family 
corresponding to a finite irreducible Coxeter group. 
These results are used to 
demonstrate that the only generalized Cartan 
matrices for which there exist finite 
edge-colored ranked posets enjoying a certain structure property 
are the Cartan matrices for the finite-dimensional complex semisimple 
Lie algebras.  
In this setting, classifications of the finite-dimensional Kac--Moody 
algebras and of the finite Coxeter and 
Weyl groups are re-derived.  
\begin{center}

\ 

\vspace*{-0.1in}

{\small \bf Keywords:}\ numbers game, generalized Cartan matrix, 
Dynkin diagram, Coxeter/Weyl group, semisimple Lie algebra, 
Kac--Moody algebra, edge-colored ranked poset  

\

\vspace*{-0.1in}

{\small \bf AMS Subject Classification:}\ 05E15 (15A48, 20F55, 17B67)

\end{center}
\end{abstract}


\noindent
{\Large \bf \IntroNum.\ \ Introduction}

\vspace{1ex} 
Dynkin diagrams are certain finite simple graphs whose edges carry certain 
information.  For some examples, see 
the connected Dynkin diagrams of ``finite type'' depicted in 
\DynkinDiagramFigure\ below. 
Many mathematical 
objects are classified by Dynkin diagrams, perhaps most famously 
the finite-dimensional complex 
semisimple Lie algebras.  For examples of other Dynkin diagram 
classifications, see 
\cite{HHSV}, \cite{PrAdv}, and \cite{PrMonthly}.  
The main results of this paper (\MainResults) are Dynkin diagram 
classifications obtained in the context of the so-called 
``numbers game.''  
These results can be viewed as classifications of 
certain combinatorial 
finiteness phenomena which are related to 
Coxeter/Weyl groups, Kac-Moody Lie algebras, and their representations. 
For example, in \S \AlgebraClassificationNum\ we recapitulate 
observations made by the second author 
in \cite{ErikssonThesis} to say how 
the Dynkin diagram classifications of 
the finite-dimensional Kac--Moody  
algebras 
(the finite-dimensional complex semisimple Lie algebras cf.\ 
\cite{Hum}, \cite{Kac}) 
and of the finite Coxeter and Weyl groups 
(see \cite{HumCoxeter}) can be re-derived in this context. 
See \ProofCorollary. For further connections, see \S \RankedPosetsNum\ 
and \S \RemarksNum\ below or \cite{DonEnumbers}. 

The numbers game is a one-player game played on a finite simple graph 
with certain real number 
weights (which we call ``amplitudes'') on its edges 
and with an initial assignment of 
real numbers 
to its nodes.  
The move a player can make 
is to ``fire'' one of the nodes with a positive number.  This move 
transforms the number at the fired node 
by changing its sign, and it also 
transforms the number at each adjacent node in a certain way 
using an amplitude   
along the incident edge.  
The player fires the nodes in some sequence of 
the player's choosing, continuing until no node has a positive 
number.  

The first main results of this paper (\FirstMainTheorems) address the 
question: 
for which such graphs 
does there exist a nontrivial initial assignment of nonnegative 
numbers such 
that the numbers game terminates in a finite number of steps?  
For graphs with negative integer amplitudes,  
we use combinatorial methods to show in \MarsAttacksTheorem\ that 
the only such connected graphs are the Dynkin 
diagrams of \DynkinDiagramFigure.   These Dynkin diagrams are in 
one-to-one correspondence  with the finite-dimensional complex simple Lie algebras. 
More generally, allowing for certain real number amplitudes we show 
in \SecondMainResult\ that the 
resulting graphs are the ``E-Coxeter graphs'' of \ECoxeterGraphFigure, 
which are related to the finite irreducible Coxeter groups. 
Our proof of this latter result 
borrows extensively from the second author's thesis \cite{ErikssonThesis} and 
applies the Perron--Frobenius theory for nonnegative real 
matrices. A different proof 
of \SecondMainResult\ 
that depends on the classification of finite 
Coxeter groups and on results 
concerning certain 
geometric representations of Coxeter groups is given in \cite{DonEnumbers}.  
All of these approaches 
build on the second author's  
work in \cite{ErikssonThesis}, 
\cite{ErikssonDiscrete}, and \cite{ErikssonEur}. 


Our third main result (\PosetTheorem) is an application of 
\MarsAttacksTheorem\ (or \SecondMainResult) to edge-colored ranked posets.  
In connection with a poset-theoretic study of the finite-dimensional 
representations 
of a given finite-dimensional complex 
semisimple Lie algebra, the first author observed in \cite{DonSupp} 
that its ``supporting 
graphs'' are finite edge-colored ranked posets possessing a certain structure 
property relative to the Cartan matrix for the Lie algebra.  
``Crystal graphs'' (see for example \cite{Stem}) and ``splitting 
posets'' (see \cite{ADLMPPW}) are other edge-colored ranked posets 
associated to semisimple Lie algebra representations which also 
possess this structure property. 
Here we begin to 
address the problem of classifying those matrices for which there are  
finite edge-colored ranked posets possessing the associated  
structure property.  In 
\PosetTheorem\ we show  
that among the ``generalized Cartan matrices'' (defined below), the 
only such matrices are those arising as Cartan matrices for 
semisimple Lie algebras.  
While our 
motivation for studying this property comes from Lie-theoretic 
considerations, 
the statement and proof of \PosetTheorem\ are purely combinatorial. 
It is hoped that, as in \cite{DonTwoColor}, this structure property 
will be part of combinatorial characterizations of families of nice 
supporting graphs or splitting posets. 

The numbers game as formulated by Mozes \cite{Mozes}  
has also been studied by Proctor 
\cite{PrEur}, \cite{PrDComplete}, Bj\"{o}rner 
\cite{Bjorner}, 
Eriksson \cite{ErikssonLinear}, \cite{ErikssonThesis}, 
\cite{ErikssonJerusalem}, \cite{ErikssonReachable}, 
\cite{ErikssonDiscrete}, \cite{ErikssonEur}, 
and Wildberger 
\cite{WildbergerAdv}, \cite{WildbergerEur}, 
\cite{WildbergerPreprint}, and Donnelly \cite{DonEnumbers}.  
Wildberger studies a dual version which 
he calls the ``mutation game.'' 
See Alon {\em et al} \cite{AKP} for a 
brief and readable treatment of the numbers game on ``unweighted'' 
cyclic graphs. Much of the numbers game discussion in Chapter 4 of 
the book \cite{BB} by Bj\"{o}rner and Brenti can be 
found in \cite{ErikssonThesis} and \cite{ErikssonDiscrete}.  
Proctor developed this process in \cite{PrEur} 
to compute Weyl group 
orbits of weights with respect to the fundamental weight basis. Here 
we use his perspective of firing 
nodes with positive, as opposed to negative, numbers. 
Mozes studied numbers games on graphs for 
which the matrix $M$ of integer amplitudes is {\em symmetrizable} 
(i.e.\ there is a nonsingular  
diagonal matrix $D$ such that $D^{-1}M$ is symmetric); 
in \cite{Mozes} he obtained ``strong convergence'' results and a 
geometric characterization of the initial positions for which the 
game terminates.  Our main results make no symmetrizable assumption.

\vspace{1ex} 

\noindent
{\Large \bf \FirstMainResultNum.\ \ 
Definitions and statement of first main result}

\vspace{1ex} 
Fix a positive integer $n$ and a totally ordered set $I_{n}$ with $n$ 
elements (usually $I_{n} := \{1<\ldots<n\}$).  
A {\em generalized Cartan matrix} (or {\em GCM}) is an $n \times n$ 
matrix $M = (M_{ij})_{i,j \in I_{n}}$  
with integer entries satisfying the requirements that each 
main diagonal matrix entry is 2, that all other matrix entries are 
nonpositive, and that if a matrix entry $M_{ij}$ is nonzero then its 
transpose entry $M_{ji}$ is also nonzero.  
Generalized Cartan matrices 
are the starting point for the study of 
Kac--Moody 
algebras: beginning with a GCM, one can 
write down a list of the defining relations for a Kac--Moody 
algebra as well as the associated Weyl group 
(see \cite{Kac}, \cite{Kumar}, or \S \AlgebraClassificationNum\ below).  
To an $n \times n$ generalized Cartan matrix 
$M = (M_{ij})_{i,j \in I_{n}}$ we associate a finite 
graph $\Gamma$ (which has undirected edges, 
no loops, and no multiple edges) 
as follows:     
The nodes $(\gamma_{i})_{i \in I_{n}}$ of $\Gamma$ are indexed 
by the set $I_{n}$, 
and   an edge is placed between nodes $\gamma_{i}$ and $\gamma_{j}$ 
if and only if $i \not= j$ 
and the matrix entries $M_{ij}$ and $M_{ji}$ are nonzero.  Call the pair 
$(\Gamma,M)$ a {\em GCM graph}. 
We consider two 
GCM graphs $(\Gamma, M = (M_{ij})_{i,j \in I_{n}})$ 
and $(\Gamma', M' = (M'_{pq})_{p,q \in I'_{n}})$ 
to be the same if under some bijection $\sigma: I_{n} \rightarrow 
I'_{n}$ we have nodes $\gamma_{i}$ and $\gamma_{j}$ in 
$\Gamma$ adjacent if and only if 
$\gamma'_{\sigma(i)}$ and $\gamma'_{\sigma(j)}$ are adjacent in $\Gamma'$ with 
$M_{ij} = M'_{\sigma(i),\sigma(j)}$. 
%
%
We depict a generic connected two-node GCM graph as 
\TwoCitiesGraphWithLabels, where $p = -M_{12}$ and $q = -M_{21}$. 
We use special names and notation to refer to 
two-node GCM 
graphs which have $p = 1$ and $q = 1$, $2$, or $3$ respectively:
\noindent
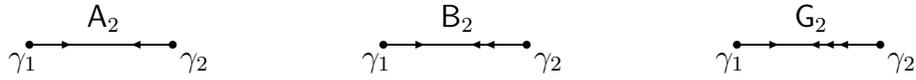
\begin{center}
\ATwoGraphWithLabels
\hspace*{1in}
\BTwoGraphWithLabels
\hspace*{1in}
\GTwoGraphWithLabels
\end{center}
When $p=1$ and $q=1$ it is convenient to use the graph 
\ATwoGraphNoEdgeLabels\ to 
represent the GCM graph $\myA_{2}$.  
A GCM graph $(\Gamma,M)$ is a {\em Dynkin diagram of finite type} if each connected 
component of $(\Gamma,M)$ 
is one of the graphs of \DynkinDiagramFigure.  Number the nodes as 
in \S 11.4 of \cite{Hum}. In these cases the GCMs  are ``Cartan'' 
matrices. 

\begin{figure}[t]
\begin{center}
\DynkinDiagramFigure: Connected Dynkin diagrams of finite type. 
\end{center}

\vspace*{-0.55in}
\begin{tabular}{cl}
$\myA_{n}$ ($n \geq 1$) & \AnMarsOK\\

$\myB_{n}$ ($n \geq 2$) & \BnMarsOK\\

$\myC_{n}$ ($n \geq 3$) & \CnMarsOK\\

$\myD_{n}$ ($n \geq 4$) & \DnMarsOK\\

$\myE_{6}$ & \ESixMarsOK\\

$\myE_{7}$ & \ESevenMarsOK\\

$\myE_{8}$ & \EEightMarsOK\\

$\myF_{4}$ & \FFourMarsOK\\

$\myG_{2}$ & \GTwoMarsOK
\end{tabular}

\vspace*{-0.15in}
\end{figure}

A {\em position} $\lambda = (\lambda_i)_{i \in 
I_{n}}$ is an assignment of real numbers to the nodes of the GCM graph 
$(\Gamma,M)$. 
The position $\lambda$ is 
{\em dominant} (respectively, {\em strongly dominant}) if 
$\lambda_{i} \geq 0$ 
(respectively $\lambda_i > 0$) for all $i \in I_{n}$; 
$\lambda$ is {\em nonzero} if at least one $\lambda_i \not= 0$. 
For $i \in I_{n}$, the {\em 
fundamental position} $\omega_i$ is the assignment of the number  
$1$ at node $\gamma_{i}$ and the number $0$ at all other nodes.  
Given a position $\lambda$ on a GCM graph $(\Gamma,M)$, to 
{\em fire} a node $\gamma_{i}$ is to change the number at each node 
$\gamma_{j}$ of $\Gamma$ by the transformation  
\[\lambda_j	 \longmapsto \lambda_j - 
M_{ij}\lambda_i,\] provided the number at node 
$\gamma_{i}$ is positive; otherwise node $\gamma_{i}$ is not allowed 
to be fired. 
Since the generalized Cartan  
matrix $M$ assigns a pair of {\em amplitudes} ($M_{ij}$ and 
$M_{ji}$) to each edge of the 
graph $\Gamma$, we sometimes refer to GCMs as 
{\em amplitude matrices}.  
The {\em numbers game} 
is the one-player game on a GCM graph $(\Gamma,M)$ in which the player 
(1) Assigns an initial position  
to the nodes of $\Gamma$; (2) Chooses a node with a positive 
number and fires the node to obtain a new position; and (3) 
Repeats step (2) for the new position if there is at least one node 
with a positive number.  
Consider now the GCM graph $\myB_{2}$.  As one can see in \BtwoMarsDemo, 
the numbers game terminates in a finite number of steps for any 
initial position and any legal sequence of node firings, 
if it is understood that the player  
will continue to fire as long as there is at least one 
node with a positive number.  In general, 
given a position $\lambda$, a {\em game sequence 
for} $\lambda$ is the (possibly empty, possibly infinite) sequence 
$(\gamma_{i_{1}}, \gamma_{i_{2}},\ldots)$, where 
$\gamma_{i_{j}}$ 
is the $j$th node that is fired in some 
numbers game with initial position $\lambda$.  
More generally, a {\em firing sequence} from some position $\lambda$ is an 
initial portion of some game sequence played from $\lambda$; the 
phrase {\em legal firing sequence} is used to emphasize that all node 
firings in the sequence are known or assumed to be possible. 
Note that a game sequence 
$(\gamma_{i_{1}}, 
\gamma_{i_{2}},\ldots,\gamma_{i_{l}})$ 
is of finite length $l$ 
(possibly with $l = 0$) if 
the number is nonpositive at each node after the $l$th firing; in 
this case the game sequence is {\em convergent} and the 
resulting position is the {\em terminal position} for the game 
sequence.  
A connected GCM graph $(\Gamma,M)$ is {\em admissible} if 
there exists a nonzero dominant initial position with a convergent 
game sequence. 

\begin{figure}[hbt]
\begin{center}
\BtwoMarsDemo: The numbers game for the GCM graph $\myB_{2}$. 

\vspace*{0.25in}
\setlength{\unitlength}{0.6in}
\begin{picture}(1.5,5.15) 
\put(0.25,5.1){\BTwoGraphForFigure}
\put(0.25,5.3){\scriptsize $a$}
\put(1.35,5.3){\scriptsize $b$}
\put(0.3,4.9){\vector(-4,-3){0.7}}
\put(1.4,4.9){\vector(4,-3){0.7}}
\put(1.75,3.8){\BTwoGraphForFigure}
\put(1.5,4.0){\scriptsize $a+2b$}
\put(2.75,4.0){\scriptsize $-b$}
\put(2.4,3.6){\vector(0,-1){0.6}}
\put(1.75,2.5){\BTwoGraphForFigure}
\put(1.45,2.7){\scriptsize $-a-2b$}
\put(2.75,2.7){\scriptsize $a+b$}
\put(2.4,2.3){\vector(0,-1){0.6}}
\put(1.75,1.2){\BTwoGraphForFigure}
\put(1.75,1.4){\scriptsize $a$}
\put(2.65,1.4){\scriptsize $-a-b$}
\put(2.1,1.0){\vector(-4,-3){0.7}}
\put(-1.25,3.8){\BTwoGraphForFigure}
\put(-1.35,4.0){\scriptsize $-a$}
\put(-0.25,4.0){\scriptsize $a+b$}
\put(-0.7,3.6){\vector(0,-1){0.6}}
\put(-1.25,2.5){\BTwoGraphForFigure}
\put(-1.5,2.7){\scriptsize $a+2b$}
\put(-0.35,2.7){\scriptsize $-a-b$}
\put(-0.7,2.3){\vector(0,-1){0.6}}
\put(-1.25,1.2){\BTwoGraphForFigure}
\put(-1.55,1.4){\scriptsize $-a-2b$}
\put(-0.15,1.4){\scriptsize $b$}
\put(-0.4,1.0){\vector(4,-3){0.7}}
\put(0.25,-0.1){\BTwoGraphForFigure}
\put(0.15,0.1){\scriptsize $-a$}
\put(1.25,0.1){\scriptsize $-b$}
\end{picture}
\end{center}
\end{figure}

The first main result of this paper 
is: 

\noindent 
{\bf \MarsAttacksTheorem}\ \ {\sl A connected GCM graph is 
admissible if and only if it 
is a connected Dynkin diagram of finite type.  
In these cases, for any given initial position every game sequence  
will converge to the same terminal position in the 
same finite number of steps.}

Our combinatorial proof of this result is given in \S \FirstProofNum. 
\SecondMainResult\ is a more general result.  

\vspace{1ex} 

\noindent
{\Large \bf \FirstProofNum.\ \ Proof of first main result}

\vspace{1ex} 
\LemmaList\ are stated and used in \cite{DonEnumbers}.  Proofs or 
references for these results are given here as well. 
Our proof of the ``only if'' direction of the first claim of 
\MarsAttacksTheorem\ uses a series of reductions that are typical in 
Dynkin diagram classification arguments.  These reductions are 
implemented using two key results (\ErikssonTheorems) 
due to the second author.  A proof of the 
remaining assertions of \MarsAttacksTheorem\ is given at the end of 
the section.  

\noindent 
{\em Proof of the ``only if'' direction of the first claim of 
\MarsAttacksTheorem:} 

\underline{Step 1: Strong convergence.}  
Following \cite{ErikssonThesis} and \cite{ErikssonEur}, we say 
the numbers game on a GCM graph $(\Gamma,M)$ is {\em strongly 
convergent} if given any initial position, any two game sequences 
either both diverge or both converge to the same terminal position in the 
same number of steps. The next result  
follows from  
Theorem 3.1 of \cite{ErikssonEur} 
(or see Theorem 3.6 of \cite{ErikssonThesis}).  

\noindent
{\bf \StrongConvergenceTheorem\ (Strong Convergence Theorem)}\ \ 
{\sl The numbers game on a connected GCM graph 
is strongly 
convergent.}  

For this part of our proof of \MarsAttacksTheorem, 
we only require the 
following weaker result, which also applies when the 
GCM graph is not connected:  

\noindent
{\bf \StrongConvergenceCorollary} \ \ {\sl For any GCM graph, 
if a game sequence for an initial position  
$\lambda$ diverges, then all game  
sequences for $\lambda$ diverge.}
 
\underline{Step 2: Comparison.}  
The next result is an immediate consequence of Theorem 4.3 of 
\cite{ErikssonThesis} or Theorem 4.5 of 
\cite{ErikssonDiscrete}.  
The proof of this result in \cite{ErikssonThesis} uses only 
combinatorial and linear algebraic methods. 
 
\noindent 
{\bf \ComparisonTheorem\ (Comparison Theorem)}\ \ 
{\sl Given a GCM graph, suppose that a game sequence 
for an initial position $\lambda = (\lambda_{i})_{i \in 
I_{n}}$ converges.  Suppose that a position $\lambda' := (\lambda'_{i})_{i 
\in I_{n}}$ has the property that $\lambda'_{i} \leq 
\lambda_{i}$ for all $i \in I_{n}$.  Then some game sequence 
for the initial position $\lambda'$ also converges.}

Let $r$ be a positive real number.  
Observe that if 
$(\gamma_{i_{1}},\ldots,\gamma_{i_{l}})$ 
is a convergent 
game sequence for an initial position $\lambda = 
(\lambda_{i})_{i \in I_{n}}$, then 
$(\gamma_{i_{1}},\ldots,\gamma_{i_{l}})$ 
is a convergent 
game sequence for the initial position $r\lambda := 
(r\lambda_{i})_{i \in I_{n}}$. 
This observation and \ComparisonTheorem\ imply the following 
result: 

\noindent 
{\bf \ComparisonCorollary}\ \ {\sl 
Let $\lambda = (\lambda_{i})_{i \in 
I_{n}}$ be a dominant initial position such that $\lambda_{j} > 0$ 
for some $j \in I_{n}$. Suppose that a game sequence for $\lambda$ 
converges.  
Then some game   
sequence for the fundamental position $\omega_{j}$ also 
converges.}

\underline{Step 3: A catalog of connected GCM graphs that are not 
admissible.}  
The following immediate consequence of \ComparisonResults\ is useful 
in the proof of \NotMarsFriendlyCatalog: 

\noindent 
{\bf \NotMarsFriendlyLemma}\ \ {\sl A GCM graph is not admissible if 
for each fundamental position there is a divergent game  
sequence.} 

\noindent 
{\bf \NotMarsFriendlyCatalog}\ \ {\sl The connected GCM graphs of 
\NotMarsFriendlyFigure\ are not admissible.}

{\em Outline of proof.} 
By \NotMarsFriendlyLemma\ it 
suffices to show that for each graph in 
\NotMarsFriendlyFigure\ and for each fundamental position, there is a 
divergent game sequence. 
In each case we 
can 
exhibit a divergent 
game sequence which is a simple pattern of node firings.  
Remarkably, 
in all cases trial and error quickly lead us to 
these patterns. 
Our goal in this 
proof 
sketch 
is not to develop any general theory for 
finding divergent game sequences for these cases, 
but rather to show that such game 
sequences can be found and presented in an elementary (though 
sometimes tedious) manner. 
For complete details see \cite{DonSupplement}. 
The $\widetilde{\myA}$, $\widetilde{\myB}$, $\widetilde{\myC}$, 
$\widetilde{\myD}$, $\widetilde{\myE}$, and $\widetilde{\myF}$ cases are 
handled using a common line of reasoning: 
A finite sequence of legal node firings is applied to a 
position whose numbers are linear expressions in an index 
variable $k$, and it is shown that the numbers for the resulting 
position are linear expressions of the same form with respect to 
the variable $k+1$.  The $\widetilde{\myG}$ cases and the families of 
small cycles are 
handled using a variation of this kind of argument: A finite sequence of 
legal node firings is applied to a generic position satisfying certain 
inequalities, and it is shown that the resulting position also 
satisfies these inequalities. 
The two paragraphs that follow  
borrow from \cite{DonSupplement} and 
demonstrate inadmissibility for some selected graphs from our list.


\begin{figure}[t]
\begin{center}
\NotMarsFriendlyFigure: Some connected GCM graphs that are not 
admissible.\\  (Figure continues on the next page.) 
\end{center}

\hspace*{0.75in}
\fbox{The $\widetilde{\myA}$ family} 
\hspace*{0.5in}
\parbox[c]{3.1in}{
\setlength{\unitlength}{0.5in}
\begin{picture}(6.1,0.7)
            \put(0,0.1){\circle*{0.075}}
            \put(1,0.1){\circle*{0.075}}
            \put(2,0.1){\circle*{0.075}}
            \put(3,0.6){\circle*{0.075}}
            \put(4,0.1){\circle*{0.075}}
            \put(5,0.1){\circle*{0.075}}
            \put(6,0.1){\circle*{0.075}}
            \put(0,0.1){\line(1,0){2}}
            \put(0,0.1){\line(6,1){3}}
            \multiput(2,0.1)(0.4,0){5}{\line(1,0){0.2}}
            \put(4,0.1){\line(1,0){2}}
            \put(6,0.1){\line(-6,1){3}}
           \end{picture}}

\vspace*{0.25in}
\hspace*{0.75in}
\fbox{The $\widetilde{\myB}$ family} 
\hspace*{0.52in}
\parbox[c]{3.1in}{
\setlength{\unitlength}{0.5in}
\begin{picture}(4.1,0.7)
            \put(0,0.6){\circle*{0.075}}
            \put(0,0.1){\circle*{0.075}}
            \put(1,0.35){\circle*{0.075}}
            \put(2,0.35){\circle*{0.075}}
            \put(4,0.35){\circle*{0.075}}
            \put(5,0.35){\circle*{0.075}}
            \put(6,0.35){\circle*{0.075}}
            \put(0,0.1){\line(4,1){1}}
            \put(0,0.6){\line(4,-1){1}}
            \put(1,0.35){\line(1,0){1}}
            \multiput(2,0.35)(0.4,0){5}{\line(1,0){0.2}}
            \put(4,0.35){\line(1,0){2}}
            \put(5.8,0.35){\vector(-1,0){0.1}}
            \put(5.2,0.35){\vector(1,0){0.1}}
            \put(5.3,0.35){\vector(1,0){0.1}}
           \end{picture}}

\vspace*{0.1in}
\hspace*{0.75in}
\hspace*{1.542in}
\parbox[c]{3.1in}{
\setlength{\unitlength}{0.5in}
\begin{picture}(4.1,0.2)
            \put(0,0.1){\circle*{0.075}}
            \put(1,0.1){\circle*{0.075}}
            \put(2,0.1){\circle*{0.075}}
            \put(4,0.1){\circle*{0.075}}
            \put(5,0.1){\circle*{0.075}}
            \put(6,0.1){\circle*{0.075}}
            \put(0,0.1){\line(1,0){2}}
            \multiput(2,0.1)(0.4,0){5}{\line(1,0){0.2}}
            \put(4,0.1){\line(1,0){2}}
            \put(0.2,0.1){\vector(1,0){0.1}}
            \put(0.8,0.1){\vector(-1,0){0.1}}
            \put(0.7,0.1){\vector(-1,0){0.1}}
            \put(5.8,0.1){\vector(-1,0){0.1}}
            \put(5.2,0.1){\vector(1,0){0.1}}
            \put(5.3,0.1){\vector(1,0){0.1}}
           \end{picture}} 

\vspace*{0.25in}
\hspace*{0.75in}
\fbox{The $\widetilde{\myC}$ family} 
\hspace*{0.525in}
\parbox[c]{3.1in}{
\setlength{\unitlength}{0.5in}
\begin{picture}(4.1,0.2)
            \put(0,0.1){\circle*{0.075}}
            \put(1,0.1){\circle*{0.075}}
            \put(2,0.1){\circle*{0.075}}
            \put(4,0.1){\circle*{0.075}}
            \put(5,0.1){\circle*{0.075}}
            \put(6,0.1){\circle*{0.075}}
            \put(0,0.1){\line(1,0){2}}
            \multiput(2,0.1)(0.4,0){5}{\line(1,0){0.2}}
            \put(4,0.1){\line(1,0){2}}
            \put(0.2,0.1){\vector(1,0){0.1}}
            \put(0.3,0.1){\vector(1,0){0.1}}
            \put(0.8,0.1){\vector(-1,0){0.1}}
            \put(5.8,0.1){\vector(-1,0){0.1}}
            \put(5.7,0.1){\vector(-1,0){0.1}}
            \put(5.2,0.1){\vector(1,0){0.1}}
           \end{picture}}  

\vspace*{0.1in}
\hspace*{0.75in}
\hspace*{1.542in}
\parbox[c]{3.1in}{
\setlength{\unitlength}{0.5in}
\begin{picture}(4.1,0.2)
            \put(0,0.1){\circle*{0.075}}
            \put(1,0.1){\circle*{0.075}}
            \put(2,0.1){\circle*{0.075}}
            \put(4,0.1){\circle*{0.075}}
            \put(5,0.1){\circle*{0.075}}
            \put(6,0.1){\circle*{0.075}}
            \put(0,0.1){\line(1,0){2}}
            \multiput(2,0.1)(0.4,0){5}{\line(1,0){0.2}}
            \put(4,0.1){\line(1,0){2}}
            \put(0.2,0.1){\vector(1,0){0.1}}
            \put(0.8,0.1){\vector(-1,0){0.1}}
            \put(0.7,0.1){\vector(-1,0){0.1}}
            \put(5.8,0.1){\vector(-1,0){0.1}}
            \put(5.7,0.1){\vector(-1,0){0.1}}
            \put(5.2,0.1){\vector(1,0){0.1}}
           \end{picture}}  

\vspace*{0.1in}
\hspace*{0.75in}
\hspace*{1.542in}
\parbox[c]{3.1in}{
\setlength{\unitlength}{0.5in}
\begin{picture}(4.1,0.7)
            \put(0,0.6){\circle*{0.075}}
            \put(0,0.1){\circle*{0.075}}
            \put(1,0.35){\circle*{0.075}}
            \put(2,0.35){\circle*{0.075}}
            \put(4,0.35){\circle*{0.075}}
            \put(5,0.35){\circle*{0.075}}
            \put(6,0.35){\circle*{0.075}}
            \put(0,0.1){\line(4,1){1}}
            \put(0,0.6){\line(4,-1){1}}
            \put(1,0.35){\line(1,0){1}}
            \multiput(2,0.35)(0.4,0){5}{\line(1,0){0.2}}
            \put(4,0.35){\line(1,0){2}}
            \put(5.2,0.35){\vector(1,0){0.1}}
            \put(5.8,0.35){\vector(-1,0){0.1}}
            \put(5.7,0.35){\vector(-1,0){0.1}}
           \end{picture}}

\vspace*{0.25in}
\hspace*{0.75in}
\fbox{The $\widetilde{\myD}$ family} 
\hspace*{0.5in}
\parbox[c]{3.1in}{
\setlength{\unitlength}{0.5in}
\begin{picture}(4.1,0.7)
            \put(0,0.6){\circle*{0.075}}
            \put(0,0.1){\circle*{0.075}}
            \put(1,0.35){\circle*{0.075}}
            \put(2,0.35){\circle*{0.075}}
            \put(4,0.35){\circle*{0.075}}
            \put(5,0.35){\circle*{0.075}}
            \put(6,0.1){\circle*{0.075}}
            \put(6,0.6){\circle*{0.075}}
            \put(0,0.1){\line(4,1){1}}
            \put(0,0.6){\line(4,-1){1}}
            \put(1,0.35){\line(1,0){1}}
            \multiput(2,0.35)(0.4,0){5}{\line(1,0){0.2}}
            \put(4,0.35){\line(1,0){1}}
            \put(5,0.35){\line(4,1){1}}
            \put(5,0.35){\line(4,-1){1}}
           \end{picture}} 

\vspace*{0.25in}
\fbox{The $\widetilde{\myE}$ family}\\ 
\hspace*{0.5in}
\parbox[c]{2.1in}{
\setlength{\unitlength}{0.5in}
\begin{picture}(4.1,1.2)
            \put(0,0.1){\circle*{0.075}}
            \put(1,0.1){\circle*{0.075}}
            \put(2,0.1){\circle*{0.075}}
            \put(2,0.6){\circle*{0.075}}
            \put(2,1.1){\circle*{0.075}}
            \put(3,0.1){\circle*{0.075}}
            \put(4,0.1){\circle*{0.075}}
            \put(0,0.1){\line(1,0){4}}
            \put(2,0.1){\line(0,1){1}}
           \end{picture}}
\hspace*{0.25in} 
\parbox[c]{3.1in}{
\setlength{\unitlength}{0.5in}
\begin{picture}(4.1,0.8)
            \put(0,0.1){\circle*{0.075}}
            \put(1,0.1){\circle*{0.075}}
            \put(2,0.1){\circle*{0.075}}
            \put(3,0.1){\circle*{0.075}}
            \put(3,0.6){\circle*{0.075}}
            \put(4,0.1){\circle*{0.075}}
            \put(5,0.1){\circle*{0.075}}
            \put(6,0.1){\circle*{0.075}}
            \put(0,0.1){\line(1,0){6}}
            \put(3,0.1){\line(0,1){0.5}}
           \end{picture}}

\vspace*{0.1in}
\hspace*{2.0in} 
\parbox[c]{3.1in}{
\setlength{\unitlength}{0.5in}
\begin{picture}(4.1,0.8)
            \put(0,0.1){\circle*{0.075}}
            \put(1,0.1){\circle*{0.075}}
            \put(2,0.1){\circle*{0.075}}
            \put(2,0.6){\circle*{0.075}}
            \put(3,0.1){\circle*{0.075}}
            \put(4,0.1){\circle*{0.075}}
            \put(5,0.1){\circle*{0.075}}
            \put(6,0.1){\circle*{0.075}}
            \put(7,0.1){\circle*{0.075}}
            \put(0,0.1){\line(1,0){7}}
            \put(2,0.1){\line(0,1){0.5}}
           \end{picture}}

\vspace*{0.25in}
\fbox{The $\widetilde{\myF}$ family}  
\hspace*{0.25in}
\parbox[c]{2.1in}{
\setlength{\unitlength}{0.5in}
\begin{picture}(4.1,0.2)
            \put(0,0.1){\circle*{0.075}}
            \put(1,0.1){\circle*{0.075}}
            \put(2,0.1){\circle*{0.075}}
            \put(3,0.1){\circle*{0.075}}
            \put(4,0.1){\circle*{0.075}}
            \put(0,0.1){\line(1,0){4}}
            \put(1.2,0.1){\vector(1,0){0.1}}
            \put(1.3,0.1){\vector(1,0){0.1}}
            \put(1.8,0.1){\vector(-1,0){0.1}}
            \end{picture}}
\hspace*{0.5in}
\parbox[c]{2.1in}{
\setlength{\unitlength}{0.5in}
\begin{picture}(4.1,0.2)
            \put(0,0.1){\circle*{0.075}}
            \put(1,0.1){\circle*{0.075}}
            \put(2,0.1){\circle*{0.075}}
            \put(3,0.1){\circle*{0.075}}
            \put(4,0.1){\circle*{0.075}}
            \put(0,0.1){\line(1,0){4}}
            \put(1.2,0.1){\vector(1,0){0.1}}
            \put(1.8,0.1){\vector(-1,0){0.1}}
            \put(1.7,0.1){\vector(-1,0){0.1}}
            \end{picture}}

\begin{center}
\underline{\hspace*{6.2in}}
\end{center}

\vspace*{-0.25in}            
\end{figure} 


In the $\widetilde{\myC}$ family, we show why 
\hspace*{0.025in}
\parbox[c]{3.1in}{
\setlength{\unitlength}{0.5in}
\begin{picture}(4.1,0.2)
            \put(0,0.1){\circle*{0.075}}
            \put(1,0.1){\circle*{0.075}}
            \put(2,0.1){\circle*{0.075}}
            \put(4,0.1){\circle*{0.075}}
            \put(5,0.1){\circle*{0.075}}
            \put(6,0.1){\circle*{0.075}}
            \put(0,0.1){\line(1,0){2}}
            \multiput(2,0.1)(0.4,0){5}{\line(1,0){0.2}}
            \put(4,0.1){\line(1,0){2}}
            \put(0.2,0.1){\vector(1,0){0.1}}
            \put(0.8,0.1){\vector(-1,0){0.1}}
            \put(0.7,0.1){\vector(-1,0){0.1}}
            \put(5.8,0.1){\vector(-1,0){0.1}}
            \put(5.7,0.1){\vector(-1,0){0.1}}
            \put(5.2,0.1){\vector(1,0){0.1}}
           \end{picture}}  
is not admissible when the graph has $n \geq 3$ nodes.  
(Since firing the middle node in the $n=3$ case 
is comparable to firing $\gamma_{2}$ in the 
$n \geq 4$ cases, then the $n=3$ case 
does not need to be considered separately here.) 
Label the nodes as $\gamma_{1}$, $\gamma_{2}$, 
$\ldots$ , $\gamma_{n-1}$, and $\gamma_{n}$ from left to right.  
For each fundamental position, we exhibit a divergent game sequence 
as a short sequence of legal 
node firings which can be repeated indefinitely. 
The fundamental position $\omega_{1} = (1,0,\ldots,0)$ is the 
$k=0$ version of the position $(2k+1,-k,0,\ldots,0)$.  From any such 
position with $k \geq 0$, 
the following sequence of node firings is easily seen to be legal: 
$(\gamma_{1}$, $\gamma_{2}$, $\ldots$ , $\gamma_{n-1}$, $\gamma_{n}$, 
$\gamma_{n-1}$, $\ldots$ , $\gamma_{3}$, 
$\gamma_{2})$.  
This sequence results in the position 
$(2(k+1)+1,-(k+1),0,\ldots,0)$. 
For $2 \leq i \leq n-1$, 
any fundamental position $\omega_{i} = (0,\ldots,0,1,0,\ldots,0)$ is the 
$k=0$ version of the position $(0,\ldots,0,2k+1,-2k,0,\ldots,0)$.  From any such 
position with $k \geq 0$, 
the following sequence of node firings is easily seen to be legal: 
$(\gamma_{i}$, $\gamma_{i-1}$, $\ldots$ , $\gamma_{2}$, 
$\gamma_{1}$, $\gamma_{2}$, 
$\ldots$ , $\gamma_{n-1}$, $\gamma_{n}$, 
$\gamma_{n-1}$, $\ldots$ , $\gamma_{i+2}$, 
$\gamma_{i+1})$.  
This sequence results in the position 
$(0,\ldots,0,2(k+1)+1,-2(k+1),0,\ldots,0)$. 
The fundamental position $\omega_{n} = (0,\ldots,0,1)$ is the 
$k=0$ version of the position $(0,\ldots,0,-2k,2k+1)$.  From any such 
position with $k \geq 0$, 
the following sequence of node firings is easily seen to be legal: 
$(\gamma_{n}$, $\gamma_{n-1}$, $\ldots$ , $\gamma_{2}$, $\gamma_{1}$, 
$\gamma_{2}$, $\ldots$ , $\gamma_{n-2}$, 
$\gamma_{n-1})$.  
This sequence results in the position 
$(0,\ldots,0,-2(k+1),2(k+1)+1)$. 


\begin{figure}[t]
\begin{center}
\NotMarsFriendlyFigure\ (continued): Some connected GCM graphs that are not 
admissible. 

\vspace*{0.25in}
\fbox{The $\widetilde{\myG}$ family}  
\hspace*{0.25in}
\parbox[c]{1.1in}{
\setlength{\unitlength}{0.5in}
\begin{picture}(2.1,0.2)
            \put(0,0.1){\circle*{0.075}}
            \put(1,0.1){\circle*{0.075}}
            \put(2,0.1){\circle*{0.075}}
            \put(0,0.1){\line(1,0){2}}
            \put(0.2,0.1){\vector(1,0){0.1}}
            \put(0.8,0.1){\vector(-1,0){0.1}}
            \put(0.7,0.1){\vector(-1,0){0.1}}
            \put(0.6,0.1){\vector(-1,0){0.1}}
            \end{picture}}
\hspace*{0.5in}
\parbox[c]{1.1in}{
\setlength{\unitlength}{0.5in}
\begin{picture}(2.1,0.2)
            \put(0,0.1){\circle*{0.075}}
            \put(1,0.1){\circle*{0.075}}
            \put(2,0.1){\circle*{0.075}}
            \put(0,0.1){\line(1,0){2}}
            \put(0.2,0.1){\vector(1,0){0.1}}
            \put(0.8,0.1){\vector(-1,0){0.1}}
            \put(0.7,0.1){\vector(-1,0){0.1}}
            \put(0.6,0.1){\vector(-1,0){0.1}}
            \put(1.2,0.1){\vector(1,0){0.1}}
            \put(1.8,0.1){\vector(-1,0){0.1}}
            \put(1.7,0.1){\vector(-1,0){0.1}}
            \end{picture}}
\hspace*{0.5in}
\parbox[c]{1.1in}{
\setlength{\unitlength}{0.5in}
\begin{picture}(2.1,0.2)
            \put(0,0.1){\circle*{0.075}}
            \put(1,0.1){\circle*{0.075}}
            \put(2,0.1){\circle*{0.075}}
            \put(0,0.1){\line(1,0){2}}
            \put(0.2,0.1){\vector(1,0){0.1}}
            \put(0.8,0.1){\vector(-1,0){0.1}}
            \put(0.7,0.1){\vector(-1,0){0.1}}
            \put(0.6,0.1){\vector(-1,0){0.1}}
            \put(1.2,0.1){\vector(1,0){0.1}}
            \put(1.8,0.1){\vector(-1,0){0.1}}
            \put(1.7,0.1){\vector(-1,0){0.1}}
            \put(1.6,0.1){\vector(-1,0){0.1}}
            \end{picture}}

\vspace*{0.1in}
\hspace*{1.28in}
\parbox[c]{1.1in}{
\setlength{\unitlength}{0.5in}
\begin{picture}(2.1,0.2)
            \put(0,0.1){\circle*{0.075}}
            \put(1,0.1){\circle*{0.075}}
            \put(2,0.1){\circle*{0.075}}
            \put(0,0.1){\line(1,0){2}}
            \put(0.2,0.1){\vector(1,0){0.1}}
            \put(0.3,0.1){\vector(1,0){0.1}}
            \put(0.4,0.1){\vector(1,0){0.1}}
            \put(0.8,0.1){\vector(-1,0){0.1}}
            \end{picture}}
\hspace*{0.5in}
\parbox[c]{1.1in}{
\setlength{\unitlength}{0.5in}
\begin{picture}(2.1,0.2)
            \put(0,0.1){\circle*{0.075}}
            \put(1,0.1){\circle*{0.075}}
            \put(2,0.1){\circle*{0.075}}
            \put(0,0.1){\line(1,0){2}}
            \put(0.2,0.1){\vector(1,0){0.1}}
            \put(0.8,0.1){\vector(-1,0){0.1}}
            \put(0.7,0.1){\vector(-1,0){0.1}}
            \put(0.6,0.1){\vector(-1,0){0.1}}
            \put(1.8,0.1){\vector(-1,0){0.1}}
            \put(1.2,0.1){\vector(1,0){0.1}}
            \put(1.3,0.1){\vector(1,0){0.1}}
            \end{picture}} 
\hspace*{0.5in}
\parbox[c]{1.1in}{
\setlength{\unitlength}{0.5in}
\begin{picture}(2.1,0.2)
            \put(0,0.1){\circle*{0.075}}
            \put(1,0.1){\circle*{0.075}}
            \put(2,0.1){\circle*{0.075}}
            \put(0,0.1){\line(1,0){2}}
            \put(0.2,0.1){\vector(1,0){0.1}}
            \put(0.8,0.1){\vector(-1,0){0.1}}
            \put(0.7,0.1){\vector(-1,0){0.1}}
            \put(0.6,0.1){\vector(-1,0){0.1}}
            \put(1.8,0.1){\vector(-1,0){0.1}}
            \put(1.2,0.1){\vector(1,0){0.1}}
            \put(1.3,0.1){\vector(1,0){0.1}}
            \put(1.4,0.1){\vector(1,0){0.1}}
            \end{picture}}

\vspace*{0.25in}
\fbox{Families of small cycles}\\
\SmallCycles  
\end{center}

\vspace*{-0.25in}            
\begin{center}
\underline{\hspace*{6.2in}}
\end{center}

\vspace*{-0.25in}            
\end{figure}


In the families of small cycles,  
we show why GCM graphs of the form \hspace*{0.1in}
\parbox[c]{0.5in}{
\setlength{\unitlength}{0.75in}
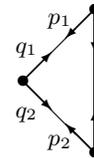
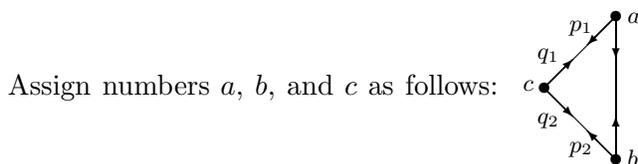
\begin{picture}(0.6,1.2)
            \put(0,0.6){\circle*{0.075}}
            \put(0.5,0.1){\circle*{0.075}}
            \put(0.5,1.1){\circle*{0.075}}
            \put(0,0.6){\line(1,1){0.5}}
            \put(0,0.6){\line(1,-1){0.5}}
            \put(0.5,0.1){\line(0,1){1}}
            \put(0.5,0.3){\vector(0,1){0.1}}
            \put(0.5,0.9){\vector(0,-1){0.1}}
            \put(0,0.6){\vector(1,1){0.2}}
            \put(-0.05,0.8){\footnotesize $q_{1}$}
            \put(0.175,1){\footnotesize $p_{1}$}
            \put(0.5,1.1){\vector(-1,-1){0.2}}
            \put(0,0.6){\vector(1,-1){0.2}}
            \put(0.5,0.1){\vector(-1,1){0.2}}
            \put(-0.05,0.35){\footnotesize $q_{2}$}
            \put(0.175,0.15){\footnotesize $p_{2}$} 
\end{picture}
}
are not admissible.  
Assign numbers $a$, $b$, 
and $c$ as follows:  
\hspace*{0.1in}
\parbox[c]{0.5in}{
\setlength{\unitlength}{0.75in}
\begin{picture}(0.6,1.2)
            \put(0,0.6){\circle*{0.075}}
            \put(0.5,0.1){\circle*{0.075}}
            \put(0.5,1.1){\circle*{0.075}}
            \put(0,0.6){\line(1,1){0.5}}
            \put(0,0.6){\line(1,-1){0.5}}
            \put(0.5,0.1){\line(0,1){1}}
            \put(0.5,0.3){\vector(0,1){0.1}}
            \put(0.5,0.9){\vector(0,-1){0.1}}
            \put(0,0.6){\vector(1,1){0.2}}
            \put(-0.05,0.8){\footnotesize $q_{1}$}
            \put(0.175,1){\footnotesize $p_{1}$}
            \put(0.5,1.1){\vector(-1,-1){0.2}}
            \put(0,0.6){\vector(1,-1){0.2}}
            \put(0.5,0.1){\vector(-1,1){0.2}}
            \put(-0.05,0.35){\footnotesize $q_{2}$}
            \put(0.175,0.15){\footnotesize $p_{2}$} 
            \put(-0.15,0.58){\footnotesize $c$}
            \put(0.58,0.05){\footnotesize $b$}
            \put(0.58,1.05){\footnotesize $a$}
\end{picture}
}
Set $\kappa := (p_{1}+p_{2}-\frac{1}{q_{2}})a + 
(p_{1}+p_{2}-\frac{1}{q_{1}})b + c$.  
Assume for now that $a \geq 0$, $b \geq 0$, $c \leq 0$, and $\kappa > 0$; 
when these inequalities hold we will say the position $(a,b,c)$ meets 
condition ({\tt *}). Under 
condition ({\tt *}) notice that $a$ and $b$ cannot both be zero. 
Begin by firing only at the two 
rightmost nodes.  When this is no longer possible, fire at the 
leftmost node.  
The resulting corresponding 
numbers are $a_{1} = 
q_{1}(\kappa + \frac{1}{q_{2}}a)$, $b_{1} = q_{2}(\kappa + 
\frac{1}{q_{1}}b)$, and $c_{1} = -\kappa-\frac{1}{q_{2}}a-\frac{1}{q_{1}}b$. In 
particular, $a_{1} > 0$, 
$b_{1} > 0$, and $c_{1} < 0$.  
Next we 
check that $\kappa_{1} := (p_{1}+p_{2}-\frac{1}{q_{2}})a_{1} + 
(p_{1}+p_{2}-\frac{1}{q_{1}})b_{1} + c_{1}$ is also positive.  
Now 
\[\kappa_{1} = Q\kappa + Q_{1}a + Q_{2}b,\]
\noindent 
where $Q = q_{1}(p_{2} - \frac{1}{q_{2}}) + q_{2}(p_{1} - 
\frac{1}{q_{1}}) + (p_{1}q_{1} + p_{2}q_{2} - 1)$, 
$Q_{1} = \frac{1}{q_{2}}[q_{1}(p_{2}-\frac{1}{q_{2}}) + 
(p_{1}q_{1} - 1)]$, and 
$Q_{2} = \frac{1}{q_{1}}[q_{2}(p_{1}-\frac{1}{q_{1}}) + 
(p_{2}q_{2} - 1)]$.  
Since each 
parenthesized quantity in our expression for $Q$ is nonnegative and 
the last of these is  
positive, then $Q > 0$.  
Similar reasoning shows that each bracketed quantity in our 
expressions for $Q_{1}$ and $Q_{2}$ is nonnegative, hence $Q_{1} \geq 0$ 
and $Q_{2} \geq 0$. 
Since $\kappa > 0$ by hypothesis, it now follows that 
$\kappa_{1} > 0$. 
Then $(a_{1},b_{1},c_{1})$ meets condition ({\tt *}), so we 
can legally 
repeat the above firing sequence from position $(a_{1},b_{1},c_{1})$ 
to obtain another position 
$(a_{2},b_{2},c_{2})$ that meets condition 
({\tt *}), etc.  
Since the fundamental 
positions $(a,b,c) = (1,0,0)$ and $(a,b,c) = (0,1,0)$ meet  
condition ({\tt *}), then we see that the indicated 
legal firing sequence can be repeated 
indefinitely from these positions.  
For the fundamental 
position $(a,b,c) = (0,0,1)$, begin by firing 
at the leftmost node to 
obtain the position $(q_{1},q_{2},-1)$. This latter position meets  
condition ({\tt *}) with $\kappa = Q$, 
and so the legal firing sequence indicated above can 
be repeated indefinitely from this position.

\underline{Step 4: Every node is fired.} The following 
is proved easily with an induction argument on the number of 
nodes.  

\noindent 
{\bf \EveryNodeFiredLemma}\ \ {\sl Suppose $(\Gamma,M)$ is connected 
with nonzero dominant position $\lambda$.  Then in any 
convergent 
game sequence for $\lambda$, every node of $\Gamma$ is fired at least 
once.}

\underline{Step 5: Subgraphs.}  
If $I'_{m}$ is a subset of the node set $I_{n}$ of a GCM graph 
$(\Gamma,M)$, 
then let $\Gamma'$ be the subgraph of $\Gamma$ with node set $I'_{m}$ and 
the induced set of edges, and let $M'$ be the corresponding 
submatrix 
of the amplitude matrix $M$; we call $(\Gamma',M')$ a {\em GCM 
subgraph} of $(\Gamma,M)$.  In light of \NewLemmaList, the 
following result amounts to an observation. 

\noindent 
{\bf \SubgraphLemma}\ \ {\sl If a connected GCM graph is 
admissible, then any 
connected GCM subgraph is also admissible.} 

\underline{Step 6: Amplitude products must be 1, 2, or 3.}

\noindent 
{\bf \TwoNodeLemma}\ \ {\sl If $\gamma_{i}$ and $\gamma_{j}$ are adjacent 
nodes in a 
connected admissible GCM graph $(\Gamma,M)$, then the product of the 
amplitudes $M_{ij}M_{ji}$ is 1, 2, or 3.  That 
is, the GCM subgraph of $(\Gamma,M)$ with nodes $\gamma_{i}$ and 
$\gamma_{j}$ 
is in this case one of} $\myA_{2}$, $\myB_{2}$, {\sl or} $\myG_{2}$.

{\em Proof.} By \SubgraphLemma\ we may restrict attention to the admissible  
GCM subgraph $(\Gamma',M')$ 
with node set $\{i,j\}$.  A nonzero dominant position with 
a convergent game sequence might not begin with 
positive numbers at 
\underline{both} nodes; nonetheless, by examining the 
proof one sees that Lemma 
3.7 of \cite{ErikssonEur} still applies to show that the product 
$M_{ij}M_{ji}$ of amplitudes in the admissible GCM 
graph $(\Gamma',M')$ is 1, 2, or 3.  

\underline{Conclusion of the ``only if'' part of the first claim of 
\MarsAttacksTheorem.} 
Putting Steps 1 through 6 together, 
we see that the only possible connected 
admissible GCM graphs are the Dynkin diagrams of finite type. 

\noindent 
{\em Proof of the remaining claims of 
\MarsAttacksTheorem:} 

Let $(\Gamma,M)$ be a connected Dynkin diagram of finite type.  
The Strong Convergence Theorem  
shows that if a game sequence for some initial position $\lambda$ converges, 
then all game sequences from $\lambda$ converge to the same 
terminal position in the 
same finite number of steps.  
Then in light of 
the Comparison Theorem, 
it  suffices to show that for 
any strongly dominant initial position on 
$(\Gamma,M)$, 
there is a convergent game sequence.  
Complete details for a case-by-case argument are given in  
\cite{DonSupplement}.  We summarize the work done there as follows.  
For the exceptional graphs 
($\myE_{6}$, $\myE_{7}$, $\myE_{8}$, $\myF_{4}$, and $\myG_{2}$) 
this can be checked by 
hand (requiring 36, 63, 120, 24, and 6 firings respectively).  
For each of the four infinite families of Dynkin diagrams of finite 
type, the proof of the next result given in \cite{DonSupplement} is 
straightforward and uses 
induction on the number of nodes. 



\noindent 
{\bf \FourFamiliesSecondLemma}\ \ {\sl For any positive integer $n$ 
(respectively, any integer $n \geq 2$, $n \geq 3$, $n \geq 4$) and for 
any strongly dominant position $(a_{1},\ldots,a_{n})$ 
on} 
$\myA_{n}$ {\sl (respectively,} $\myB_{n}${\sl ,} $\myC_{n}${\sl ,} 
$\myD_{n}${\sl ), one can obtain the 
terminal position $(-a_{n},\ldots,-a_{2},-a_{1})$ 
(respectively 
$(-a_{1},-a_{2},\ldots,-a_{n})$, 
$(-a_{1},-a_{2},\ldots,-a_{n})$, 
$(-a_{1},-a_{2},\ldots,
-a_{n-2},-b_{n-1},-b_{n})$ where $b_{n-1} := a_{n-1}$ and $b_{n} := 
a_{n}$ when $n$ is even  
and where $b_{n-1} := a_{n}$ and $b_{n} := a_{n-1}$ 
when $n$ is odd) by a sequence of $\frac{n(n+1)}{2}$ (resp.\ 
$n^{2}$, $n^{2}$, $n(n-1)$) node firings.} 

This completes the proof of \MarsAttacksTheorem.\hfill\QED

\vspace{1ex} 

\noindent
{\Large \bf \SecondMainResultNum.\ \ A more general admissibility 
result}

\vspace{1ex} 
In this section we prove  
a generalization (\SecondMainResult) 
of \MarsAttacksTheorem: Allowing for certain real 
number amplitudes, we show that the connected 
admissible graphs are precisely those depicted in 
\ECoxeterGraphFigure.  
The proof of \SecondMainResult\  
uses reasoning from Ch.\ 5 and 6 of 
\cite{ErikssonThesis}. 
The key is the application of the 
Perron--Frobenius theory of nonnegative real matrices to a matrix 
closely related to the amplitude matrix for the graph. 

Readers familiar with Coxeter groups will notice that the families of 
graphs of 
\ECoxeterGraphFigure\ correspond to the 
finite irreducible Coxeter groups.  
However, the proof of \SecondMainResult\ requires no Coxeter group theory.  
We will see in \S \AlgebraClassificationNum\ 
how the classification of finite irreducible 
Coxeter groups can be deduced from \SecondMainResult.  

Now we play the numbers game in the following more 
general environment. 
An {\em E-generalized Cartan matrix} or {\em E-GCM} is  
a real $n \times n$ matrix $M = (M_{ij})_{i,j \in I_{n}}$  
satisfying the requirements that each 
main diagonal matrix entry is 2, that all other matrix entries are 
nonpositive, that if a matrix entry $M_{ij}$ is nonzero then its 
transpose entry $M_{ji}$ is also nonzero, and that if 
$M_{ij}M_{ji}$ is nonzero then $M_{ij}M_{ji} \geq 4$ or 
$M_{ij}M_{ji} = 4\cos^{2}(\pi/k_{ij})$ for some integer $k_{ij} \geq 
3$.  
We apply the language of GCM graphs in this E-GCM setting, so an {\em 
E-GCM graph} is a pair $(\Gamma,M)$ where $M$ is an E-GCM, etc.  
For the remainder of this section, $(\Gamma,M)$ denotes an E-GCM graph. 
Examining the proofs and references given in \S \FirstProofNum, 
one sees that 
\EGCMLemmaList\ as well as the Strong Convergence and Comparison 
Theorems
stated for GCM graphs also hold for E-GCM graphs.  
Any particular two-node E-GCM graph is depicted in the same way as a 
two-node GCM graph.  We use \EGCMGraphCirclem 
for the collection of all two-node E-GCM 
graphs for which $M_{12}M_{21} = pq = 4\cos^{2}(\pi/m)$ 
for an integer $m > 3$; we use $m = \infty$ if $M_{12}M_{21} = pq \geq 
4$.  When $m = 3$ (i.e.\ $pq = 1$), 
we use an unlabelled edge $\!\!$\ATwoGraphNoEdgeLabels\ $\!\!$.  
An {\em E-Coxeter graph} will be any E-GCM graph whose 
connected components come from one of the collections of 
\ECoxeterGraphFigure. 

The peculiar constraints on products of transpose pairs of E-GCM  
entries are precisely those required in order to guarantee strong 
convergence for these ``E-games'' 
(see \cite{ErikssonEur}).  
These constraints also afford a precise connection 
between E-games and certain geometric representations of Coxeter 
groups.  This connection was developed by Vinberg \cite{Vinberg} and 
the second author   
\cite{ErikssonThesis}, \cite{ErikssonDiscrete}, and it has also 
been studied in \cite{BB}, \cite{DonEnumbers}, and \cite{ProctorCoxeter}.  

\begin{figure}[t]
\begin{center}
\ECoxeterGraphFigure: Families of connected E-Coxeter graphs. 
\vspace*{0.04in}

{\footnotesize (For adjacent nodes, the notation$\!\!$ 
\CircleIntegerm$\!\!$ means 
that the amplitude product on the edge is $4\cos^{2}(\pi/m)$;\\ for an 
unlabelled edge take $m=3$.)}
\end{center}

\vspace*{-0.45in}
\begin{tabular}{cl}
$\mathcal{A}_{n}$ ($n \geq 1$) & \AnEGraph\\

$\mathcal{B}_{n}$ ($n \geq 3$) & \BnEGraph\\

$\mathcal{D}_{n}$ ($n \geq 4$) & \DnEGraph\\

$\mathcal{E}_{6}$ & \ESixEGraph\\

$\mathcal{E}_{7}$ & \ESevenEGraph\\

$\mathcal{E}_{8}$ & \EEightEGraph\\

$\mathcal{F}_{4}$ & \FFourEGraph\\

$\mathcal{H}_{3}$ & \HThreeEGraph\\

$\mathcal{H}_{4}$ & \HFourEGraph\\

$\mathcal{I}_{2}^{(m)}$ ($4 \leq m < \infty$) & \ITwoEGraph
\end{tabular}
\end{figure}


The main result of this section is: 

\noindent
{\bf \SecondMainResult}\ \ {\sl A connected E-GCM graph  is 
admissible if and only if it 
is a connected E-Coxeter graph.  
In these cases, for any given initial position every game sequence  
will converge to the same terminal position in the 
same finite number of steps.}  

The proof is at the end of the section and requires the following 
linear algebra set-up.  A matrix (or vector) of real numbers is {\em 
nonnegative} if all its entries are nonnegative.  For real matrices 
$X = (X_{ij})$ and $Y = (Y_{ij})$, 
say $X \geq Y$ if 
$X_{ij} \geq Y_{ij}$ for all $i, j$.  With respect to this 
partial ordering, we let $X > Y$ mean $X \geq 
Y$ but $X \not= Y$.  A square matrix $X = (X_{ij})_{i,j \in I_{n}}$ 
is {\em indecomposable} if there is no permutation of $I_{n}$ such 
that we get a block matrix $\left(\begin{array}{cc}X^{(1)} & X^{(2)}\\ 
O & X^{(3)}\end{array}\right)$ where $X^{(1)}$ and $X^{(3)}$ are 
square and $O$ denotes a zero matrix.  A {\em principal submatrix} is 
a submatrix obtained by deleting some rows and the corresponding 
columns.  Let $\varrho(X)$ denote the spectral radius, i.e.\ the 
modulus of the largest eigenvalue, of a square matrix $X$. The 
well-known results comprising the following  theorem 
can be found for example in \cite{Minc}. 

\noindent 
{\bf Theorem (Perron--Frobenius)}\ \ {\sl A nonnegative square 
indecomposable matrix $X$ has one eigenvector, unique up to  
multiplication with a scalar, with all elements positive.  Its 
corresponding eigenvalue is $\varrho(X) > 0$, and it is simple.  No 
other eigenvector is nonnegative.  If $Y$ is a principal submatrix 
of $X$, then $\varrho(X) > \varrho(Y)$.  If $X > Y$, then again $\varrho(X) > 
\varrho(Y)$.}

For $J \subseteq I_{n}$ and 
$m = |J|$, let $E(J)$ denote 
the $m \times m$ identity matrix $(\delta_{ij})_{i,j \in J}$.  
Let $E := E(I_{n})$. 
Let $\Gamma$ be a finite simple graph with node set 
$\{\gamma_{x}\}_{x \in I_{n}}$.  
Let $A$ be a nonnegative matrix such that 
$A_{ij} > 0$ if and only if $A_{ji} > 0$ if and only if 
$\gamma_{i}$ and $\gamma_{j}$ are 
adjacent nodes in $\Gamma$.  Let $A^{\mbox{\scriptsize sym}}$ denote the 
symmetric nonnegative 
matrix for which $A^{\mbox{\scriptsize sym}}_{ij}A^{\mbox{\scriptsize 
sym}}_{ji} = A_{ij}A_{ji}$ for all $i,j$.  

\noindent
{\bf \EGCMLinAlgLemma}\ \ {\sl Keep the notation of the preceding paragraph.  (1) 
$\Gamma$ is connected if and only if $A$ is indecomposable. (2) If $\Gamma$ 
is acyclic, then $A$ has the 
same characteristic polynomial as} $A^{\mbox{\scriptsize sym}}${\sl , and 
all eigenvalues of $A$ are real. (3) Let} $B := \mytranspose{A}-2E$.  
{\sl If $\Gamma$ is connected and 
acyclic and if 
$\varrho(A) < 2$, then 
there is a positive definite diagonal matrix $D$ such that $DB^{-1}$ is 
symmetric and positive definite.} 

{\em Proof.} For (1), $A$ is indecomposable if and only if there is 
no partitioning $I_{n} = J_{1} \cup J_{2}$ such that $A_{ij} = 0$ 
when $i \in J_{1}$ and $j \in J_{2}$.  That is, $A$ is indecomposable 
iff $\Gamma$ is connected.  For (2), induct on the number of nodes. 
Without loss of generality assume 
node $\gamma_{1}$ has at most one adjacent node in $\Gamma$.  
Let $\Gamma'$ be the 
subgraph obtained by removing node $\gamma_{1}$.  Let $A'$ be the 
corresponding principal submatrix of $A$, and set $J' := I_{n} 
\setminus \{1\}$.  If a node $\gamma_{2}$ is 
adjacent to $\gamma_{1}$, let 
$\Gamma''$ be the subgraph obtained by removing $\gamma_{1}$ and 
$\gamma_{2}$, and let $A''$ be the corresponding principal submatrix.  
In any case, set $J'' := I_{n} 
\setminus \{1,2\}$. Let $x$ be an indeterminate.  It is easy to check that 
\[\det(xE-A) = x\det(xE(J')-A') - A_{12}A_{21}\det(xE(J'')-A'').\]
Assuming the result holds for such graphs on fewer than $n$ nodes, 
it follows from this computation that $\det(xE-A) = 
\det(xE-A^{\mbox{\scriptsize sym}})$.  
Since a symmetric real matrix has real eigenvalues, the result follows. 
For (3), note that all eigenvalues 
for $B$ are real and positive.  By \S 1.5 of \cite{Kumar}, 
$\mytranspose{B}$ is symmetrizable and  
there is a positive definite diagonal matrix $D$ such that 
$D^{-1}\mytranspose{B}$ is symmetric.  (See remarks following 
Definition 1.5.1 as well as Exercise 1.5.1 in \cite{Kumar}.)  Then 
$(D^{-1}\mytranspose{B})^{-1}$ is symmetric and positive definite, 
hence $DB^{-1}$ is symmetric and positive definite.\hfill\QED

We apply these linear algebra facts to E-games as follows. 
Think of a position $\lambda = 
(\lambda_{i})_{i \in I_{n}}$ on $(\Gamma,M)$ 
as an $n \times 1$ column vector.  
Firing node $\gamma_{i}$ from $\lambda$ results in the position 
$\lambda^{\mbox{\scriptsize new}} = \lambda - 
\lambda_{i}\mytranspose{M}\omega_{i}$.  
For any connected E-GCM subgraph $(\Gamma',M')$ with nodes indexed by 
the set $J := \{x \in I_{n}\}_{\gamma_{x}\in\Gamma'}$, set 
$A := A_{\Gamma',M'} := 2E(J)-M'$.  
Then $A$ is a nonnegative real matrix.  
Connectedness of $\Gamma'$ implies that $A$ is indecomposable.  
By Perron--Frobenius, there is a vector $\nu := \nu_{\Gamma',M'} = 
(\nu_{i})_{i \in I_{n}}$ 
(unique up to a positive scalar multiple) with $\nu_{i} > 0$ for all 
$i$ and $A\nu = 
\varrho(A)\nu$.  A {\em looping game} is a nonempty 
legal sequence of node 
firings from some position $\lambda$ that returns to position 
$\lambda$. 

\noindent 
{\bf \EGCMVectorLemma}\ \ {\sl 
Suppose $(\Gamma,M)$ is connected.  Let $A := A_{\Gamma,M}$ and 
$\nu := \nu_{\Gamma,M}$. (1) 
Suppose that $A$ has largest eigenvalue $\varrho(A) \geq 2$,   
that from a position $\lambda$ some node is legally fired to 
obtain the position} $\lambda^{\mbox{\scriptsize new}}${\sl , and 
that}  
$\mytranspose{\nu}\lambda > 0${\sl . Then} 
$\mytranspose{\nu}\lambda^{\mbox{\scriptsize new}} > 0${\sl , and   
some node may be legally fired from position} 
$\lambda^{\mbox{\scriptsize new}}${\sl . (2) If there is a looping game 
on $(\Gamma,M)$, then $\varrho(A) = 2$.} 

{\em Proof.} For (1) we have 
$\mytranspose{\nu}\lambda^{\mbox{\scriptsize new}} = 
\mytranspose{\nu}\lambda - \lambda_{i}\mytranspose{(M\nu)}\omega_{i} = 
\mytranspose{\nu}\lambda + (\varrho(A)-2)\lambda_{i}\nu_{i}$, 
assuming $\gamma_{i}$ is the fired node of the hypothesis. 
Since $\mytranspose{\nu}\lambda > 0$, $\lambda_{i} > 0$, 
$\nu_{i} > 0$, 
and $\varrho(A) - 2 \geq 
0$, then $\mytranspose{\nu}\lambda^{\mbox{\scriptsize new}} > 0$.  
That some some number 
$\lambda^{\mbox{\scriptsize new}}_{j}$ is positive 
follows from the facts that 
$\mytranspose{\nu}\lambda^{\mbox{\scriptsize new}} > 0$ and $\nu_{k} > 
0$ for all $k$.  For (2), let $(\gamma_{i_{1}}, \ldots, 
\gamma_{i_{j}})$ be the sequence of nodes fired for a looping game 
played from position $\lambda^{(1)}$.  Let $\lambda^{(1)}$, $\lambda^{(2)}$, 
$\ldots$, $\lambda^{(j)}$, $\lambda^{(j+1)}$ be the sequence of 
positions for the looping game, with $\lambda^{(j+1)} = 
\lambda^{(1)}$.  As in the proof of (1), we compute that 
$\mytranspose{\nu}\lambda^{(j+1)} = 
\mytranspose{\nu}\lambda^{(1)} +  
(\varrho(A)-2)\sum_{k=1}^{j}\lambda^{(k)}_{i_{k}}\nu_{i_{k}}$. 
Then 
$(\varrho(A)-2)\sum_{k=1}^{j}\lambda^{(k)}_{i_{k}}\nu_{i_{k}} = 
0$. Since $\sum_{k=1}^{j}\lambda^{(k)}_{i_{k}}\nu_{i_{k}} > 0$, then 
$\varrho(A) = 2$.\hfill\QED 

We apply part (1) of the preceding 
lemma to obtain the following criterion for inadmissibility. 

\noindent 
{\bf \EGCMInadmissibleProp}\ \ {\sl Suppose $(\Gamma,M)$ is connected, and 
suppose $(\Gamma',M')$ is a connected E-GCM subgraph with $J := \{x 
\in I_{n}\}_{\gamma_{x}\in\Gamma'}$.  Suppose $A' := A_{\Gamma',M'}$ 
has largest eigenvalue $\varrho(A') \geq 2$. 
Then $(\Gamma,M)$ is not admissible.} 

{\em Proof.} We prove the contrapositive. 
Suppose $(\Gamma,M)$ is admissible.  Then for some 
nonzero dominant position $\lambda$, there is a convergent game 
sequence. The Strong Convergence Theorem guarantees that all game 
sequences played from $\lambda$ 
are convergent.  In particular, we may choose a convergent 
game sequence that starts by firing only nodes indexed by the set 
$I_{n} \setminus J$.  When there are no positive numbers on nodes 
for $I_{n} \setminus J$, we then have a dominant position $\lambda'$ 
on $(\Gamma',M')$.  Since every node must be fired in a convergent 
game sequence played from $\lambda$ (\EveryNodeFiredLemma), then 
$\lambda'$ must be nonzero.  
Using Perron--Frobenius, take a vector $\nu' = 
(\nu_{j}')_{j \in J}$ such that $\nu_{j}' > 0$ for all $j$ and 
$A'\nu' = \varrho(A')\nu'$.  Since $\lambda'$ is nonzero and dominant 
for $(\Gamma',M')$ and since $\nu_{j}' > 0$ for all $j$, 
it follows that $\mytranspose{(\nu')}\lambda' > 0$.  
If $\varrho(A') \geq 2$, then \EGCMVectorLemma\  
shows we will have a divergent game sequence on $(\Gamma',M')$ played 
from $\lambda'$, and 
hence a divergent game sequence on $(\Gamma,M)$ 
from $\lambda$.  So it must be the 
case that $\varrho(A') < 2$.\hfill\QED 



\noindent
{\bf \EGCMAdmissibleProp}\ \ {\sl Suppose $(\Gamma,M)$ is a connected 
E-GCM graph for which $\varrho(A_{\Gamma,M}) < 2$. Let $\lambda$ be 
any position on $(\Gamma,M)$.  Then there is a convergent game 
sequence played from $\lambda$.} 

{\em Proof.} 
Set $A := A_{\Gamma,M}$  
and $B := -\mytranspose{M}$.  
From \EGCMLinAlgLemma.3 it 
follows that there 
is a positive definite diagonal matrix $D$ such that $DB^{-1}$ is 
symmetric and positive definite.  
Suppose that a position $\mu^{\mbox{\scriptsize new}}$ is obtained 
from some position $\mu$ by legally firing $\gamma_{i}$.  We claim that 
$\mytranspose{(\mu^{\mbox{\scriptsize new}})}DB^{-1}\mu^{\mbox{\scriptsize 
new}} = \mytranspose{\mu}DB^{-1}\mu$.  To see this, make the 
substitution  
$\mu^{\mbox{\scriptsize new}} = \mu+\mu_{i}B\omega_{i}$ on the 
left-hand side.  Then the claim is equivalent to the statement that 
$\mu_{i}\mytranspose{\mu}D\omega_{i} 
+ \mu_{i}\mytranspose{\omega_{i}}\mytranspose{B}DB^{-1}\mu  
+ \mu_{i}^{2}\mytranspose{\omega_{i}}\mytranspose{B}D\omega_{i} = 0$. 
Now $\mu_{i}\mytranspose{\mu}D\omega_{i}$ is just 
$\mu_{i}^{2}D_{ii}$.  Since $DB^{-1}$ is symmetric, then 
$\mu_{i}\mytranspose{\omega_{i}}\mytranspose{B}DB^{-1}\mu$ is also 
$\mu_{i}^{2}D_{ii}$.  Finally, 
$\mu_{i}^{2}\mytranspose{\omega_{i}}\mytranspose{B}D\omega_{i} = 
\mu_{i}^{2}D_{ii}B_{ii} = -2\mu_{i}^{2}D_{ii}$. 

Topologize the space of positions by identifying it with 
$\mathbb{R}^{n}$.  
Now let $\lambda$ be any position.  
Let $\mathfrak{P}(\lambda)$ be the set of all positions that are 
reachable from position $\lambda$ by legal firing sequences.  Note 
that $\mathfrak{P}(\lambda)$ is a discrete subset of the space of all 
positions.  The 
preceding paragraph shows that $\mathfrak{P}(\lambda) \subseteq 
\{\mbox{positions } \mu\, |\, \mytranspose{\mu}DB^{-1}\mu = 
\mytranspose{\lambda}DB^{-1}\lambda\}$.  The latter set is compact 
since $DB^{-1}$ is positive definite.  
Since $\mathfrak{P}(\lambda)$ is a discrete subset of 
a compact set, then $\mathfrak{P}(\lambda)$ is finite.  By 
\EGCMVectorLemma.2  
there are no looping games on $(\Gamma,M)$.  Since there are no 
repeated positions in E-game play from $\lambda$ and since the set of 
positions reachable from $\lambda$ is finite, it follows that there 
must be a convergent game sequence from $\lambda$.\hfill\QED 

\noindent 
{\bf \EGCMInadmissibleList}\ \ {\sl Let $(\Gamma,M)$ be an E-GCM graph 
from one of the families of \EGCMInadmissibleFigure. If $(\Gamma,M)$ is 
from one 
of the $\tilde{\mathcal{H}}_{3}$ or $\tilde{\mathcal{H}}_{3}$
families, then $\varrho(A_{\Gamma,M}) > 2$.  
If $(\Gamma,M)$ is in 
the $\tilde{\mathcal{A}}_{n}$ 
family $(n \geq 1)$, then $\varrho(A_{\Gamma,M}) \geq 2$.  
If $(\Gamma,M)$ is from any 
other 
family of \EGCMInadmissibleFigure, then $\varrho(A_{\Gamma,M}) = 2$.} 

\begin{figure}[ht]
\begin{center}
\EGCMInadmissibleFigure: Some connected E-GCM graphs for 
\EGCMInadmissibleList
\vspace*{0.025in}

{\footnotesize (Notation: An E-GCM graph in family 
$\tilde{\mathcal{X}}_{n}$ below has $n+1$ nodes.\\ The numbers 
assigned to the 
nodes of some of these graphs pertain to the proof of \EGCMInadmissibleList.)}
\end{center}

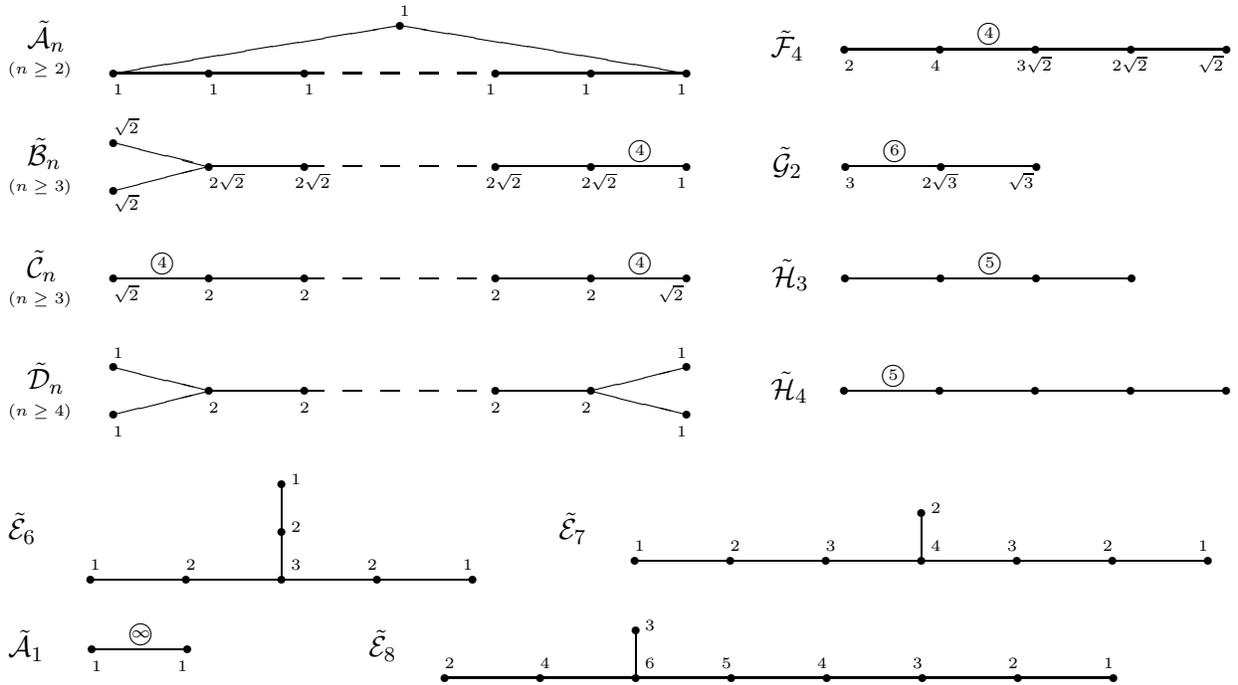
\parbox{0.6in}{\hspace*{0.1in}$\tilde{\mathcal{A}}_{n}$ 

\vspace*{-0.05in}
{\tiny $(n \geq 2)$}}
\hspace*{-0.15in}
\parbox[c]{3.1in}{
\setlength{\unitlength}{0.5in}
\begin{picture}(6.1,0.7)
            \put(0,0.1){\circle*{0.075}}
            \put(1,0.1){\circle*{0.075}}
            \put(2,0.1){\circle*{0.075}}
            \put(3,0.6){\circle*{0.075}}
            \put(4,0.1){\circle*{0.075}}
            \put(5,0.1){\circle*{0.075}}
            \put(6,0.1){\circle*{0.075}}
            \put(0,0.1){\line(1,0){2}}
            \put(0,0.1){\line(6,1){3}}
            \multiput(2,0.1)(0.4,0){5}{\line(1,0){0.2}}
            \put(4,0.1){\line(1,0){2}}
            \put(6,0.1){\line(-6,1){3}}
            \put(0,-0.125){\tiny $1$}
            \put(1,-0.125){\tiny $1$}
            \put(2,-0.125){\tiny $1$}
            \put(3.9,-0.125){\tiny $1$}
            \put(4.9,-0.125){\tiny $1$}
            \put(5.9,-0.125){\tiny $1$}
            \put(3,0.7){\tiny $1$}
           \end{picture}}
\hspace*{0.24in}
$\tilde{\mathcal{F}}_{4}$
\hspace*{0.11in}
\parbox[c]{2.1in}{
\setlength{\unitlength}{0.5in}
\begin{picture}(4.1,0.2)
            \put(0,0.1){\circle*{0.075}}
            \put(1,0.1){\circle*{0.075}}
            \put(2,0.1){\circle*{0.075}}
            \put(3,0.1){\circle*{0.075}}
            \put(4,0.1){\circle*{0.075}}
            \put(0,0.1){\line(1,0){4}}
            \put(1.2,0.15){\CircleInteger{4}}
            \put(0,-0.125){\tiny $2$}
            \put(0.9,-0.125){\tiny $4$}
            \put(1.8,-0.125){\tiny $3\sqrt{2}$}
            \put(2.8,-0.125){\tiny $2\sqrt{2}$}
            \put(3.7,-0.125){\tiny $\sqrt{2}$}
            \end{picture}}

\vspace*{0.25in}
\parbox{0.6in}{\hspace*{0.1in}$\tilde{\mathcal{B}}_{n}$

\vspace*{-0.05in}
{\tiny $(n \geq 3)$}}
\hspace*{-0.15in}
\parbox[c]{3.1in}{
\setlength{\unitlength}{0.5in}
\begin{picture}(4.1,0.7)
            \put(0,0.6){\circle*{0.075}}
            \put(0,0.1){\circle*{0.075}}
            \put(1,0.35){\circle*{0.075}}
            \put(2,0.35){\circle*{0.075}}
            \put(4,0.35){\circle*{0.075}}
            \put(5,0.35){\circle*{0.075}}
            \put(6,0.35){\circle*{0.075}}
            \put(0,0.1){\line(4,1){1}}
            \put(0,0.6){\line(4,-1){1}}
            \put(1,0.35){\line(1,0){1}}
            \multiput(2,0.35)(0.4,0){5}{\line(1,0){0.2}}
            \put(4,0.35){\line(1,0){2}}
            \put(5.2,0.4){\CircleInteger{4}}
            \put(0,-0.1){\tiny $\sqrt{2}$}
            \put(0,0.7){\tiny $\sqrt{2}$}
            \put(1,0.125){\tiny $2\sqrt{2}$}
            \put(1.9,0.125){\tiny $2\sqrt{2}$}
            \put(3.9,0.125){\tiny $2\sqrt{2}$}
            \put(4.9,0.125){\tiny $2\sqrt{2}$}
            \put(5.9,0.125){\tiny $1$}
           \end{picture}}
\hspace*{0.25in}
$\tilde{\mathcal{G}}_{2}$
\hspace*{0.125in}
\parbox[c]{1.1in}{
\setlength{\unitlength}{0.5in}
\begin{picture}(2.1,0.2)
            \put(0,0.1){\circle*{0.075}}
            \put(1,0.1){\circle*{0.075}}
            \put(2,0.1){\circle*{0.075}}
            \put(0,0.1){\line(1,0){2}}
            \put(0.2,0.15){\CircleInteger{6}}
            \put(0,-0.125){\tiny $3$}
            \put(0.8,-0.125){\tiny $2\sqrt{3}$}
            \put(1.7,-0.125){\tiny $\sqrt{3}$}
            \end{picture}}

\vspace*{0.25in}
\parbox{0.6in}{\hspace*{0.1in}$\tilde{\mathcal{C}}_{n}$

\vspace*{-0.05in}
{\tiny $(n \geq 3)$}}
\hspace*{-0.15in}
\parbox[c]{3.1in}{
\setlength{\unitlength}{0.5in}
\begin{picture}(4.1,0.2)
            \put(0,0.1){\circle*{0.075}}
            \put(1,0.1){\circle*{0.075}}
            \put(2,0.1){\circle*{0.075}}
            \put(4,0.1){\circle*{0.075}}
            \put(5,0.1){\circle*{0.075}}
            \put(6,0.1){\circle*{0.075}}
            \put(0,0.1){\line(1,0){2}}
            \multiput(2,0.1)(0.4,0){5}{\line(1,0){0.2}}
            \put(4,0.1){\line(1,0){2}}
            \put(0.2,0.15){\CircleInteger{4}}
            \put(5.2,0.15){\CircleInteger{4}}
            \put(0,-0.125){\tiny $\sqrt{2}$}            
            \put(0.95,-0.125){\tiny $2$}            
            \put(1.95,-0.125){\tiny $2$}            
            \put(3.95,-0.125){\tiny $2$}            
            \put(4.95,-0.125){\tiny $2$}            
            \put(5.7,-0.125){\tiny $\sqrt{2}$}            
           \end{picture}}  
\hspace*{0.245in}
$\tilde{\mathcal{H}}_{3}$
\hspace*{0.09in}
\parbox[c]{2.1in}{
\setlength{\unitlength}{0.5in}
\begin{picture}(3.1,0.2)
            \put(0,0.1){\circle*{0.075}}
            \put(1,0.1){\circle*{0.075}}
            \put(2,0.1){\circle*{0.075}}
            \put(3,0.1){\circle*{0.075}}
            \put(0,0.1){\line(1,0){3}}
            \put(1.2,0.15){\CircleInteger{5}}
            \end{picture}}

\vspace*{0.25in}
\parbox{0.6in}{\hspace*{0.1in}$\tilde{\mathcal{D}}_{n}$ 

\vspace*{-0.05in}
{\tiny $(n \geq 4)$}}
\hspace*{-0.15in}
\parbox[c]{3.1in}{
\setlength{\unitlength}{0.5in}
\begin{picture}(4.1,0.7)
            \put(0,0.6){\circle*{0.075}}
            \put(0,0.1){\circle*{0.075}}
            \put(1,0.35){\circle*{0.075}}
            \put(2,0.35){\circle*{0.075}}
            \put(4,0.35){\circle*{0.075}}
            \put(5,0.35){\circle*{0.075}}
            \put(6,0.1){\circle*{0.075}}
            \put(6,0.6){\circle*{0.075}}
            \put(0,0.1){\line(4,1){1}}
            \put(0,0.6){\line(4,-1){1}}
            \put(1,0.35){\line(1,0){1}}
            \multiput(2,0.35)(0.4,0){5}{\line(1,0){0.2}}
            \put(4,0.35){\line(1,0){1}}
            \put(5,0.35){\line(4,1){1}}
            \put(5,0.35){\line(4,-1){1}}
            \put(0,-0.125){\tiny $1$}
            \put(0,0.7){\tiny $1$}
            \put(1,0.125){\tiny $2$}
            \put(1.95,0.125){\tiny $2$}
            \put(3.95,0.125){\tiny $2$}
            \put(4.9,0.125){\tiny $2$}
            \put(5.9,-0.125){\tiny $1$}
            \put(5.9,0.7){\tiny $1$}
           \end{picture}} 
\hspace*{0.24in}
$\tilde{\mathcal{H}}_{4}$
\hspace*{0.09in}
\parbox[c]{2.1in}{
\setlength{\unitlength}{0.5in}
\begin{picture}(4.1,0.2)
            \put(0,0.1){\circle*{0.075}}
            \put(1,0.1){\circle*{0.075}}
            \put(2,0.1){\circle*{0.075}}
            \put(3,0.1){\circle*{0.075}}
            \put(4,0.1){\circle*{0.075}}
            \put(0,0.1){\line(1,0){4}}
            \put(0.2,0.15){\CircleInteger{5}}
            \end{picture}}

\vspace*{0.25in}
$\tilde{\mathcal{E}}_{6}$
\hspace*{0.185in}
\parbox[c]{2.1in}{
\setlength{\unitlength}{0.5in}
\begin{picture}(4.1,1.2)
            \put(0,0.1){\circle*{0.075}}
            \put(1,0.1){\circle*{0.075}}
            \put(2,0.1){\circle*{0.075}}
            \put(2,0.6){\circle*{0.075}}
            \put(2,1.1){\circle*{0.075}}
            \put(3,0.1){\circle*{0.075}}
            \put(4,0.1){\circle*{0.075}}
            \put(0,0.1){\line(1,0){4}}
            \put(2,0.1){\line(0,1){1}}
            \put(0,0.2){\tiny $1$}
            \put(1,0.2){\tiny $2$}
            \put(2.1,1.1){\tiny $1$}
            \put(2.1,0.6){\tiny $2$}
            \put(2.1,0.2){\tiny $3$}
            \put(2.9,0.2){\tiny $2$}
            \put(3.9,0.2){\tiny $1$}
           \end{picture}}
\hspace*{0.25in} 
$\tilde{\mathcal{E}}_{7}$
\hspace*{0.15in} 
\parbox[c]{3.1in}{
\setlength{\unitlength}{0.5in}
\begin{picture}(4.1,0.8)
            \put(0,0.1){\circle*{0.075}}
            \put(1,0.1){\circle*{0.075}}
            \put(2,0.1){\circle*{0.075}}
            \put(3,0.1){\circle*{0.075}}
            \put(3,0.6){\circle*{0.075}}
            \put(4,0.1){\circle*{0.075}}
            \put(5,0.1){\circle*{0.075}}
            \put(6,0.1){\circle*{0.075}}
            \put(0,0.1){\line(1,0){6}}
            \put(3,0.1){\line(0,1){0.5}}
            \put(0,0.2){\tiny $1$}
            \put(1,0.2){\tiny $2$}
            \put(2,0.2){\tiny $3$}
            \put(3.1,0.6){\tiny $2$}
            \put(3.1,0.2){\tiny $4$}
            \put(3.9,0.2){\tiny $3$}
            \put(4.9,0.2){\tiny $2$}
            \put(5.9,0.2){\tiny $1$}
           \end{picture}}

\vspace*{0.1in}
$\tilde{\mathcal{A}}_{1}$
\hspace*{0.15in}
\parbox[c]{1.1in}{
\setlength{\unitlength}{0.5in}
\begin{picture}(1.1,0.2)
            \put(0,0.1){\circle*{0.075}}
            \put(1,0.1){\circle*{0.075}}
            \put(0,0.1){\line(1,0){1}}
            \put(0.2,0.15){\CircleInfty}
            \put(0,-0.125){\tiny $1$}
            \put(0.9,-0.125){\tiny $1$}
            \end{picture}}
\hspace*{0.25in} 
$\tilde{\mathcal{E}}_{8}$
\hspace*{0.15in} 
\parbox[c]{3.1in}{
\setlength{\unitlength}{0.5in}
\begin{picture}(4.1,0.8)
            \put(0,0.1){\circle*{0.075}}
            \put(1,0.1){\circle*{0.075}}
            \put(2,0.1){\circle*{0.075}}
            \put(2,0.6){\circle*{0.075}}
            \put(3,0.1){\circle*{0.075}}
            \put(4,0.1){\circle*{0.075}}
            \put(5,0.1){\circle*{0.075}}
            \put(6,0.1){\circle*{0.075}}
            \put(7,0.1){\circle*{0.075}}
            \put(0,0.1){\line(1,0){7}}
            \put(2,0.1){\line(0,1){0.5}}
            \put(0,0.2){\tiny $2$}
            \put(1,0.2){\tiny $4$}
            \put(2.1,0.6){\tiny $3$}
            \put(2.1,0.2){\tiny $6$}
            \put(2.9,0.2){\tiny $5$}
            \put(3.9,0.2){\tiny $4$}
            \put(4.9,0.2){\tiny $3$}
            \put(5.9,0.2){\tiny $2$}
            \put(6.9,0.2){\tiny $1$}
           \end{picture}}
\end{figure}

{\em Proof.} Consider an E-GCM graph $(\Gamma,M)$ from one of the 
$\tilde{\mathcal{B}}_{n}$, $\tilde{\mathcal{C}}_{n}$, 
$\tilde{\mathcal{D}}_{n}$, $\tilde{\mathcal{E}}_{6}$, 
$\tilde{\mathcal{E}}_{7}$, $\tilde{\mathcal{E}}_{8}$,  
$\tilde{\mathcal{F}}_{4}$, or $\tilde{\mathcal{G}}_{2}$ families of 
\EGCMInadmissibleFigure.  Take $M$ to be symmetric. 
Regard the numbers on the nodes as coordinates for a vector 
$\nu$.  
One can confirm in each case that 
$A_{\Gamma,M}\nu = 2\nu$, and hence that 
$\varrho(A_{\Gamma,M}) = 2$ and $\nu = \nu_{\Gamma,M}$.  
By \EGCMLinAlgLemma.2, it follows that 
$\varrho(A_{\Gamma',M'}) = 2$ for any E-GCM graph $(\Gamma',M')$ 
in one of these families. 
Now take $(\Gamma,M)$ in the 
$\tilde{\mathcal{A}}_{1}$ family with $M$ symmetric and $M_{12}M_{21} = 
4$.  With $\nu$ as prescribed in  \EGCMInadmissibleFigure, we get 
$A_{\Gamma,M}\nu = 2\nu$, and hence  
$\varrho(A_{\Gamma,M}) = 2$ and $\nu = \nu_{\Gamma,M}$. 
By \EGCMLinAlgLemma.2, it follows that 
$\varrho(A_{\Gamma',M'}) = 2$ for any E-GCM graph $(\Gamma',M')$ 
in the $\tilde{\mathcal{A}}_{1}$ family for which $M'_{12}M'_{21} = 4$. 
If $(\Gamma'',M'')$ is any E-GCM graph from the $\tilde{\mathcal{A}}_{1}$ 
family, then we can find a graph $(\Gamma',M')$ from this same family 
such that $M'_{12} \geq M''_{12}$, $M'_{21} \geq M''_{21}$, and 
$M'_{12}M'_{21} = 4$.  Then $A_{\Gamma',M'} \leq A_{\Gamma'',M''}$, 
and hence by Perron--Frobenius $\varrho(A_{\Gamma'',M''}) \geq 2$.  
We used a computer algebra system to calculate the eigenvalues for 
$A_{\Gamma,M}$ when $(\Gamma,M)$ is 
in one of the $\tilde{\mathcal{H}}_{3}$ or $\tilde{\mathcal{H}}_{4}$ 
families with $M$ 
symmetric.  In 
these cases one finds that $\varrho(A_{\Gamma,M}) > 2$.  It 
follows from \EGCMLinAlgLemma.2 that  
$\varrho(A_{\Gamma,M}) > 2$ for any E-GCM 
graph $(\Gamma,M)$ in the $\tilde{\mathcal{H}}_{3}$ or 
$\tilde{\mathcal{H}}_{4}$ 
families.  Now consider an E-GCM graph $(\Gamma,M)$ in the 
$\tilde{\mathcal{A}}_{n}$ family ($n\geq{2}$).  
Assume the nodes are numbered consecutively around the 
cycle.  Set $\Pi := (-1)^{n}M_{12}M_{23}\cdots{M_{n-1,n}}M_{n,1}$, so 
$\Pi > 0$.  
One can show by a 
computation that the 
characteristic polynomial of $A := A_{\Gamma,M}$ is the 
characteristic polynomial of  
$A^{\mbox{\scriptsize sym}}$ plus the constant $2-\Pi-\frac{1}{\Pi}$.  
We have $A^{\mbox{\scriptsize sym}}\nu = 2\nu$ for the vector $\nu$ 
identified in \EGCMInadmissibleFigure, so $2$ is an eigenvalue for 
$A^{\mbox{\scriptsize sym}}$.  
Since $2-\Pi-\frac{1}{\Pi}$ is nonpositive, 
it follows that $A$ has a real eigenvalue no less than $2$.\hfill\QED 

\noindent 
{\bf Remark}\ \ In \cite{ErikssonThesis}, the second author 
classifies the 
``E-loopers,'' i.e.\ those connected E-GCM graphs with looping games.  
It is shown that an E-looper must be one of: 
an E-GCM graph in one of the $\tilde{\mathcal{B}}_{n}$, 
$\tilde{\mathcal{C}}_{n}$
$\tilde{\mathcal{D}}_{n}$, $\tilde{\mathcal{E}}_{6}$, 
$\tilde{\mathcal{E}}_{7}$, $\tilde{\mathcal{E}}_{8}$, 
$\tilde{\mathcal{F}}_{4}$, or  $\tilde{\mathcal{G}}_{2}$ families; 
an E-GCM graph $(\Gamma,M)$ in the $\tilde{\mathcal{A}}_{1}$ family 
such that $-M_{12} = -M_{21} = 2$; or an E-GCM graph 
$(\Gamma,M)$ in the $\tilde{\mathcal{A}}_{n}$ family ($n \geq 2$) 
such that $(-1)^{n}M_{12}M_{23}{\cdots}M_{n-1,n}M_{n,1} = 1$, where the nodes 
of $\Gamma$ are are numbered consecutively around the 
cycle.  

\noindent
{\bf \EGCMAdmissibleLemma}\ \ 
{\sl If $(\Gamma,M)$ is a connected E-Coxeter graph, then 
$\varrho(A_{\Gamma,M}) < 2$.} 

{\em Proof.} Let $(\Gamma,M)$ be the E-GCM graph with symmetric $M$ 
from the $\mathcal{I}_{2}(m)$ family.  Let 
$(\tilde{\Gamma},\tilde{M})$ be the E-GCM graph from the 
$\tilde{\mathcal{A}}_{1}$ family with symmetric 
$\tilde{M}$ satisfying $\tilde{M}_{12}\tilde{M}_{21} = 4$. Check that 
$\varrho(A_{\tilde{\Gamma},\tilde{M}}) = 2$ (see e.g.\ the proof of  
\EGCMInadmissibleList).  
Since $A_{\Gamma,M} < A_{\tilde{\Gamma},\tilde{M}}$, then by  
Perron--Frobenius 
$\varrho(A_{\Gamma,M}) < 2 = \varrho(A_{\tilde{\Gamma},\tilde{M}})$. 
By \EGCMLinAlgLemma.2, $\varrho(A_{\Gamma,M}) < 2$ for any E-GCM 
graph $(\Gamma,M)$ in the $\mathcal{I}_{2}(m)$ family.  Use similar 
reasoning to see that 
other E-Coxeter graphs from \ECoxeterGraphFigure\ have dominant eigenvalue 
less than 2.  In particular, let $(\Gamma,M)$ be the E-Coxeter graph from family 
$\mathcal{X}_{n} \in \{\mathcal{A}_{n}, \mathcal{B}_{n}, 
\mathcal{D}_{n}, \mathcal{E}_{6}, \mathcal{E}_{7}, \mathcal{E}_{8}, 
\mathcal{F}_{4}\}$ such that $M$ is symmetric.    
Let $(\tilde{\Gamma},\tilde{M})$ be the respective 
E-Coxeter graph from family 
$\tilde{\mathcal{X}}_{n} \in \{\tilde{\mathcal{A}}_{n}, 
\tilde{\mathcal{B}}_{n}, 
\tilde{\mathcal{D}}_{n}, \tilde{\mathcal{E}}_{6}, 
\tilde{\mathcal{E}}_{7}, \tilde{\mathcal{E}}_{8}, 
\tilde{\mathcal{F}}_{4}\}$ such that $\tilde{M}$ is symmetric.  By 
Perron--Frobenius and the fact that $A_{\Gamma,M}$ is a principal 
submatrix of  
$A_{\tilde{\Gamma},\tilde{M}}$ we get $\varrho(A_{\Gamma,M}) < 2 = 
\varrho(A_{\tilde{\Gamma},\tilde{M}})$.  By \EGCMLinAlgLemma.2, 
$\varrho(A_{\Gamma,M}) < 2$ for any E-GCM 
graph $(\Gamma,M)$ in the $\mathcal{X}_{n}$ family. 
We used a computer algebra system to calculate the eigenvalues for 
$A_{\Gamma,M}$ when $(\Gamma,M)$ is 
in one of the $\mathcal{H}_{3}$ or $\mathcal{H}_{4}$ families with $M$ 
symmetric.  In 
these cases one finds that $\varrho(A_{\Gamma,M}) < 2$ as well.  It 
follows from \EGCMLinAlgLemma.2 that  
$\varrho(A_{\Gamma,M}) < 2$ for any E-Coxeter  
graph $(\Gamma,M)$ in the $\mathcal{H}_{3}$ or $\mathcal{H}_{4}$ 
families.\hfill\QED 

{\em Proof of \SecondMainResult.} 
\EGCMLemmaListForProof\ demonstrate  
admissibility for any connected E-Coxeter 
graph.  In such a case it follows from the Strong Convergence and Comparison 
Theorems that for any given initial position any two game sequences  
will converge to the same terminal position in the 
same finite number of steps.  Now 
suppose a connected E-GCM graph $(\Gamma,M)$ is admissible.  First, we show 
that $(\Gamma,M)$ has no cycles.  If an E-GCM subgraph $(\Gamma',M')$ 
is a cycle, then we can let $M''$ be an E-GCM for which $M''_{ij} \geq 
M'_{ij}$ for all appropriate $i,j$ and such that $(\Gamma'',M'')$ is in the 
$\tilde{\mathcal{A}}_{k}$ family.  Then $A_{\Gamma'',M''} < 
A_{\Gamma',M'}$, and by Perron--Frobenius, $2 \leq \varrho(A_{\Gamma'',M''}) < 
\varrho(A_{\Gamma',M'})$.  But then \EGCMInadmissibleProp\ implies 
that $(\Gamma,M)$ is not admissible.  Hence $(\Gamma,M)$ has no cycles. 
Similar reasoning using inadmissibility of E-GCM graphs in 
the $\tilde{\mathcal{B}}_{k}$ family 
shows that $(\Gamma,M)$ cannot have a node with more than two 
neighbors if $M_{ij}M_{ji} > 1$ for some pair of adjacent nodes 
$\gamma_{i}$ and $\gamma_{j}$, i.e.\ $m_{ij} > 3$.  
Assume $(\Gamma,M)$ has a pair of adjacent nodes $\gamma_{i}$ and 
$\gamma_{j}$ with $M_{ij}M_{ji} > 1$.  Inadmissibility of E-GCM 
graphs in the $\tilde{\mathcal{A}}_{1}$ family shows that for all such 
pairs of nodes $M_{ij}M_{ji} < 4$.  
Now using inadmissibility of E-GCM graphs in 
the $\tilde{\mathcal{C}}_{k}$ family 
shows that $(\Gamma,M)$ can have at most one such pair of adjacent 
nodes.  Using inadmissibility of E-GCM graphs in 
the $\tilde{\mathcal{F}}_{4}$ family 
shows that if neither $\gamma_{i}$ nor $\gamma_{j}$ is a leaf, then 
$(\Gamma,M)$ must be in the $\mathcal{F}_{4}$ family.   
If at least one of $\gamma_{i}$ or $\gamma_{j}$ is a leaf, 
then inadmissibility of E-GCM graphs in 
the $\tilde{\mathcal{G}}_{2}$, $\tilde{\mathcal{H}}_{3}$, and 
$\tilde{\mathcal{H}}_{4}$ families 
shows that 
$(\Gamma,M)$ must be in the $\mathcal{B}_{n}$, $\mathcal{C}_{n}$, 
$\mathcal{H}_{3}$,  $\mathcal{H}_{4}$,  or $\mathcal{I}_{2}(m)$   
family.   
Now assume that for each pair of adjacent nodes $\gamma_{i}$ and 
$\gamma_{j}$ in $(\Gamma,M)$ we have $M_{ij}M_{ji} = 1$, i.e.\ 
$m_{ij} = 3$.  
Using inadmissibility of E-GCM graphs in 
the $\tilde{\mathcal{D}}_{k}$ family 
shows that $(\Gamma,M)$ can have at most one node with more than two 
neighbors and that such a node can have at most three neighbors.  
Using inadmissibility of E-GCM graphs in 
the $\tilde{\mathcal{E}}_{6}$, $\tilde{\mathcal{E}}_{7}$, and 
$\tilde{\mathcal{E}}_{8}$ families 
we conclude if $(\Gamma,M)$ has a node with three neighbors such that 
at most one of them is a leaf in the tree, then $(\Gamma,M)$ must be in 
one of the $\mathcal{E}_{6}$, $\mathcal{E}_{7}$, and 
$\mathcal{E}_{8}$ families.  Otherwise, $(\Gamma,M)$ is in the 
$\mathcal{A}_{n}$ or $\mathcal{D}_{n}$ family.\hfill\QED

\vspace{1ex} 

\noindent
{\Large \bf \RankedPosetsNum.\ \ A structure property for 
edge-colored ranked posets} 

\vspace{1ex} 
The so-called structure property we study in this section is motivated 
by a certain Lie-theoretic phenomenon.  
``Crystal graphs'' (see \cite{Stem}), ``supporting graphs'' (see 
\cite{DonSupp}), and ``splitting posets'' (see \cite{ADLMPPW}) 
are edge-colored 
directed graphs which encode certain information about the 
finite-dimensional representations of a given finite-dimensional complex 
simple Lie algebra $\mathfrak{g}$.  
Ignoring edge colors, any such graph is the 
Hasse diagram for a ranked partially ordered set (defined below).  
The number $n$ of edge colors is the rank of $\mathfrak{g}$, 
i.e.\ the dimension of a Cartan subalgebra.  
The edges and edge colors for such a poset 
determine a ``weight'' rule, that is, an assignment of an integer 
$n$-tuple to each vertex of the Hasse diagram.  The weight rule has 
the following property: the difference of 
the weights for two vertices in such a graph is the $i$th row of the Cartan matrix 
for $\mathfrak{g}$ if the vertices form an edge  
whose color corresponds to $i$. 
The results of this section begin to address  
the question: If an edge-colored ranked poset possesses such a 
property relative to some matrix $M$, what can be said about $M$? 
 
The set-up of this paragraph follows Section 2 of \cite{ADLMPPW}. 
Identify a partially ordered set $R$ with its {\em Hasse diagram}, 
that is, the 
directed graph whose edges depict the {\em covering relations} for the 
poset: for elements $\selt$ and $\telt$ in $R$ the directed edge 
$\selt \rightarrow \telt$ means that $\selt < \telt$ and 
if $\selt \leq \xelt \leq \telt$ 
then $\selt = \xelt$ or $\xelt = \telt$. 
The edge set $\mathcal{E}(R)$ is the set of all covering relations 
in $R$. Given an $n$-element set $I_{n}$, a function 
$\ecolor_{R}\, :\, \mathcal{E}(R) \longrightarrow I_{n}$ is an {\em 
edge coloring function}, in which 
case we say the edges of $R$ are colored by the set $I_{n}$.  
The notation $\selt \myarrow{i} \telt$ means that $\ecolor_{R}(\selt 
\rightarrow \telt) = i$. 
If $J$ is a subset of $I_{n}$,
remove all edges from $R$ whose colors are not in $J$; connected
components of the resulting edge-colored poset are called
{\em J-components} of $R$. 
For any $\xelt$ in $R$, we let 
$\mathbf{comp}_{J}(\xelt)$ denote the $J$-component of $R$ 
containing $\xelt$.  
When $R$ is finite we say it 
is {\em ranked} if there exists a surjective function $\rho : R 
\longrightarrow \{0,1,\ldots,l\}$ such that $\rho(\selt) + 1 = 
\rho(\telt)$ whenever $\selt \rightarrow \telt$; in this case the 
number $l$ is the {\em length} of $R$, $\rho$ is a {\em rank function}, 
and $\rho(\xelt)$ is the {\em rank} of any element $\xelt$ in $R$.  
With respect to an edge coloring function 
$\ecolor_{R}\, :\, \mathcal{E}(R) \longrightarrow I_{n}$ on our 
finite ranked poset $R$, we let $\rho_{i}(\xelt)$ denote the rank of 
an element $\xelt$ in $R$ within its $i$-component 
$\mathbf{comp}_{i}(\xelt)$ and we let $l_{i}(\xelt)$ denote the length 
of $\mathbf{comp}_{i}(\xelt)$.  
For any $\xelt \in R$, let 
$wt_{R}(\xelt)$ be the $n$-tuple $(\, m_{i}(\xelt)\, )_{i \in 
I_{n}}$, where $m_{i}(\xelt) := 2\rho_{i}(\xelt) - l_{i}(\xelt)$ for 
each $i \in I_{n}$.  
Now let $M = (M_{i,j})_{i,j \in I_{n}}$ be a GCM, and for each $i \in 
I_{n}$ let $\alpha_{i}$ denote the $i$th row of $M$: $\alpha_{i} = 
(M_{i,j})_{j \in I_{n}}$.  An $M$-{\em structure poset} $(R,\ecolor_{R})$ 
is a finite ranked poset $R$ together with an edge coloring 
function $\ecolor_{R}\, :\, \mathcal{E}(R) \longrightarrow I_{n}$ 
satisfying the following $M$-{\em structure 
property}: $wt_{R}(\selt) + \alpha_{i} = wt_{R}(\telt)$ whenever 
$\selt \myarrow{i} \telt$ in $R$. Note the finiteness requirement 
of the definition. 

\noindent 
{\bf \PosetTheoremConnectedGraph}\ \ 
{\sl Suppose $M = (M_{i,j})_{i,j \in I_{n}}$ is a 
generalized Cartan matrix with connected GCM graph 
$(\Gamma,M)$, and 
suppose there is an $M$-structure poset $(R,\ecolor_{R})$ 
with at least one edge.  
Then $(\Gamma,M)$ is a 
connected Dynkin 
diagram of finite type, and $\ecolor_{R}$ is surjective.} 

{\em Proof.} Choose a vertex $\telt_{0}$ for which 
$\lambda^{(0)} := wt_{R}(\telt_{0})$ is dominant.  
(For example, take $\telt_{0}$ to 
be any element of highest rank in $R$.)  
Since $R$ has at least one edge, $\lambda^{(0)}$ is nonzero. 
Let $(\gamma_{i_{1}}, 
\gamma_{i_{2}}, \ldots)$ be any game sequence played from initial 
position $\lambda^{(0)}$ on $(\Gamma,M)$.  
For each $p \geq 1$, $\lambda^{(p)}$ is the position in the sequence 
just after node $\gamma_{i_{p}}$ is fired.   
Next, we define by induction a special sequence of elements from $R$.  
For any $p 
\geq 1$, suppose we have a sequence 
$\telt_{0}, \telt_{1}, 
\ldots, \telt_{p-1}$ for which $wt_{R}(\telt_{q}) = \lambda^{(q)}$ 
and $\rho(\telt_{q}) < \rho(\telt_{q-1})$ for all $1 \leq q \leq p-1$.  
We wish to show that we can extend this sequence by an element 
$\telt_{p}$ so that $wt_{R}(\telt_{p}) = \lambda^{(p)}$ and 
$\rho(\telt_{p}) < \rho(\telt_{p-1})$.   
Take $\telt_{p}$ to be any element of $\mathbf{comp}_{i_{p}}(\telt_{p-1})$ 
for which 
$\rho_{i_{p}}(\telt_{p}) = l_{i_{p}}(\telt_{p-1}) - 
\rho_{i_{p}}(\telt_{p-1})$. 
Since firing node $\gamma_{i_{p}}$ in the 
given numbers game is legal from position $\lambda^{(p-1)}$, 
then $\lambda^{(p-1)}_{i_{p}} > 0$.  But 
$\lambda^{(p-1)}_{i_{p}} = 2\rho_{i_{p}}(\telt_{p-1}) - 
l_{i_{p}}(\telt_{p-1})$.  So, $\rho_{i_{p}}(\telt_{p}) = 
l_{i_{p}}(\telt_{p-1}) - \rho_{i_{p}}(\telt_{p-1}) < 
\rho_{i_{p}}(\telt_{p-1})$.  It follows that $\rho(\telt_{p}) < 
\rho(\telt_{p-1})$.  
Since $R$ satisfies the $M$-structure condition, then 
$wt_{R}(\telt_{p}) = wt_{R}(\telt_{p-1}) - 
\lambda^{(p-1)}_{i_{p}}\alpha_{i_{p}}$.  But $\lambda^{(p-1)} = 
wt_{R}(\telt_{p-1})$ and $\lambda^{(p)} = 
\lambda^{(p-1)} - \lambda^{(p-1)}_{i_{p}}\alpha_{i_{p}}$.  In other 
words, $wt_{R}(\telt_{p}) = \lambda^{(p)}$.  So we have extended our 
sequence as desired. 
But since $R$ is finite, any such sequence must also be finite.  Hence 
the game sequence $(\gamma_{i_{1}}, \gamma_{i_{2}}, \ldots)$ is 
convergent.  Then by \MarsAttacksTheorem, 
$(\Gamma,M)$ must be a Dynkin diagram of finite type.  Since every 
node must be fired in a convergent game sequence for  
a numbers game played on a connected GCM graph (\EveryNodeFiredLemma), 
it follows that 
$\ecolor_{R}$ is surjective.\hfill\QED

Given an $M$-structure poset $(R,\ecolor_{R})$ for a generalized 
Cartan matrix $M = (M_{i,j})_{i,j \in I_{n}}$,  we 
say $\ecolor_{R}$ is {\em sufficiently surjective} if in each 
connected component of the GCM graph $(\Gamma,M)$ there is some node 
$\gamma_{j}$ such that $j \in \ecolor_{R}(\mathcal{E}(R))$.  

\noindent 
{\bf \PosetTheorem}\ \ {\sl Let $(\Gamma,M 
= (M_{i,j})_{i,j \in I_{n}})$ be a GCM graph.  Suppose 
there is an $M$-structure poset $(R,\ecolor_{R})$ 
with sufficiently surjective edge 
coloring function $\ecolor_{R}$. Then 
$(\Gamma,M)$ is a Dynkin 
diagram of finite type, and $\ecolor_{R}$ is surjective.} 

{\em Proof.} Pick a connected component $(\Gamma',M')$ of 
$(\Gamma,M)$,  
and let $J := \{x \in I_{n}\}_{\gamma_{x} \in \Gamma'}$. 
Now pick a 
$J$-component $\mathcal{C}$ of $R$ such that $\mathcal{C}$ contains at 
least one edge whose color is from $J$. 
\PosetTheoremConnectedGraph\ implies 
that $(\Gamma',M')$ is a connected Dynkin diagram of finite type 
and that for every 
color in $J$ there is an edge in $\mathcal{C}$ having that color.  
Applying this reasoning to each connected component of $(\Gamma,M)$, 
we see that $(\Gamma,M)$ is a Dynkin diagram of finite type 
and that $\ecolor_{R}$ is 
surjective.\hfill\QED

\noindent 
{\bf Remark}\ \ The existence of $M$-structure posets meeting the 
hypotheses of 
\PosetTheorem\ can be 
seen as follows.  Given a Dynkin diagram $(\Gamma,M)$ of finite type, 
let $\mathfrak{g}$ be the finite-dimensional complex semisimple Lie 
algebra with Cartan matrix $M$.  
Let $V$ be any 
finite-dimensional irreducible 
representation of $\mathfrak{g}$ whose highest weight has the 
following property: it is a 
nonnegative integer linear combination of fundamental weights such 
that for each connected component $(\Gamma',M')$ 
of $(\Gamma,M)$ there is a 
fundamental weight appearing nontrivially in the linear combination 
whose corresponding fundamental 
position is in $(\Gamma',M')$.  Then any supporting graph 
for $V$ is a connected $M$-structure poset satisfying 
the hypotheses of \PosetTheorem, cf.\ Lemmas 3.1 and 3.2 
of \cite{DonSupp}. 
From \S 2 of \cite{Stem}, one can similarly see that any ``admissible system'' for $V$ 
(e.g.\ a crystal graph) is a connected $M$-structure poset meeting 
the hypotheses of \PosetTheorem.

\vspace{1ex} 

\noindent
{\Large \bf \AlgebraClassificationNum.\ \ Classifications of finite-dimensional 
Kac--Moody algebras and finite Coxeter and Weyl groups}

\vspace{1ex} 
Here we re-derive some well-known Dynkin diagram 
classification results using an argument from 
\S 8.4 of \cite{ErikssonThesis}. 
Since the 
classifications of the finite-dimensional Kac--Moody 
algebras and of 
the finite Coxeter and Weyl groups are not used 
in the proofs of \FirstMainTheorems,  
we can use these theorems to obtain these 
classification results.  
This is recorded below as \ProofCorollary. 
These classifications are obtained 
in \cite{Kac} and \cite{HumCoxeter} respectively by carefully 
studying properties of the generalized Cartan matrix 
(or a closely related 
matrix).   
The key observation is that a divergent game sequence requires 
an infinite group: any finite firing sequence corresponds to an 
element of the same finite length in a corresponding Coxeter group.  
From this reasoning we obtain for free 
the classification of the finite-dimensional Kac--Moody algebras. 

Given an E-GCM graph $(\Gamma,M)$, 
define the associated Coxeter group 
$W = W(\Gamma,M)$ to be the Coxeter group with identity 
denoted $\varepsilon$, 
generators $\{s_{i}\}_{i \in I_{n}}$, and defining relations $s_{i}^{2} = 
\varepsilon$ for $i \in I_{n}$ and 
$(s_{i}s_{j})^{m_{ij}} = \varepsilon$ for all $i \not= j$, where the 
$m_{ij}$ are determined as follows: 

\vspace*{-0.1in}
\[m_{ij} = \left\{\begin{array}{cl}
k_{ij} & \mbox{\hspace*{0.25in} if 
$M_{ij}M_{ji} = 4\cos^{2}(\pi/k_{ij})$ for some integer $k_{ij} \geq 2$}\\ 
\infty & \mbox{\hspace*{0.25in} if 
$M_{ij}M_{ji} \geq 4$} 
\end{array}\right.\]  
(Conventionally, $m_{ij} = \infty$ means there is no relation between 
generators 
$s_{i}$ and $s_{j}$.) 
One can think of the E-GCM graph as a refinement of the information 
from the Coxeter graph for the associated Coxeter group.    
Observe that any Coxeter group on a finite set of generators is isomorphic 
to the Coxeter group associated to some E-GCM graph. 
The Coxeter group $W$ is {\em irreducible} if $\Gamma$ is connected.   
If the graph $\Gamma$ has connected components $\Gamma_{1}, 
\ldots, \Gamma_{k}$ with corresponding amplitude matrices 
$M_{1}, \ldots, M_{k}$, then  $W(\Gamma,M) \approx 
W(\Gamma_{1},M_{1}) \times \cdots 
\times W(\Gamma_{k},M_{k})$. 
Let ${\myl}$ denote the length function for $W$. 
An expression $s_{i_{p}}{\cdots}s_{i_{2}}s_{i_{1}}$ for an element of 
$W$ is {\em reduced} 
if $\myl(s_{i_{p}}{\cdots}s_{i_{2}}s_{i_{1}}) = p$. An empty product 
in $W$ is taken as $\varepsilon$. 
For a firing sequence $(\gamma_{i_{1}}, 
\gamma_{i_{2}}, \ldots, \gamma_{i_{p}})$ from some initial position 
on $(\Gamma,M)$, 
the corresponding element of $W$ is taken to be  
$s_{i_{p}}\cdots{s}_{i_{2}}s_{i_{1}}$. 
The next result follows from Propositions 4.1 and 4.2 of 
\cite{ErikssonDiscrete} and is the key step in the proof of 
\ProofCorollary. 


\noindent 
{\bf \ErikssonWordProposition}\ \ {\sl (1) If 
$(\gamma_{i_{1}}, 
\gamma_{i_{2}}, \ldots, \gamma_{i_{p}})$ is a legal sequence of node 
firings in a numbers game played  
from some initial position on an E-GCM graph $(\Gamma,M)$, 
then $s_{i_{p}}\cdots{s}_{i_{2}}s_{i_{1}}$ is a reduced expression for 
the corresponding element of $W(\Gamma,M)$. (2) If 
$s_{i_{p}}\cdots{s}_{i_{2}}s_{i_{1}}$ is a reduced expression for an 
element of $W(\Gamma,M)$, then $(\gamma_{i_{1}}, 
\gamma_{i_{2}}, \ldots, \gamma_{i_{p}})$ is a legal sequence of node 
firings in a numbers game played  
from any strongly dominant position on E-GCM graph $(\Gamma,M)$.} 

The Lie algebra that is constructed next 
does not depend on the specific choices made.  
The definitions we use here basically follow \cite{Kumar} (but see 
also \cite{Kac}).  
Given a GCM graph $(\Gamma,M)$ with $n$ nodes, choose a complex vector space 
$\mathfrak{h}$ of dimension $n + \mathrm{corank}(M)$.  Choose $n$ 
linearly independent vectors $\{\beta_{i}^{\vee}\}_{1 \leq i \leq n}$ 
in $\mathfrak{h}$, and find $n$ linearly independent functionals 
$\{\beta_{i}\}_{1 \leq i \leq n}$ in $\mathfrak{h}^{*}$ satisfying 
$\beta_{j}(\beta_{i}^{\vee}) = M_{ij}$.  
The {\sl Kac--Moody 
algebra} $\mathfrak{g} = 
\mathfrak{g}(\Gamma,M)$ is the Lie 
algebra over $\mathbb{C}$ generated by the set $\mathfrak{h} \cup 
\{x_{i},y_{i}\}_{i \in I_{n}}$ with relations 
$[\mathfrak{h},\mathfrak{h}] = 0$; 
$[h,x_{i}] = \beta_{i}(h)x_{i}$ and $[h,y_{i}] = 
-\beta_{i}(h)y_{i}$ for all $h \in \mathfrak{h}$ and $i \in I_{n}$; 
$[x_{i},y_{j}] = \delta_{i,j}\beta_{i}^{\vee}$ for all $i, j \in I_{n}$; 
$(\mathrm{ad} x_{i})^{1-M_{ji}}(x_{j}) = 0$ for $i \not= j$;  
and 
$(\mathrm{ad} y_{i})^{1-M_{ji}}(y_{j}) = 0$ for $i \not= j$, 
where $(\mathrm{ad} z)^{k}(w) = [z,[z,\cdots,[z,w]\cdots]]$.  
If the graph $\Gamma$ has connected components $\Gamma_{1}, 
\ldots, \Gamma_{k}$ with corresponding amplitude matrices 
$M_{1}, \ldots, M_{k}$, then 
$\mathfrak{g}(\Gamma,M) \approx \mathfrak{g}(\Gamma_{1},M_{1}) \oplus \cdots 
\oplus \mathfrak{g}(\Gamma_{k},M_{k})$. 
Note that if $M_{ij}M_{ji} < 4$, then $M_{ij}M_{ji} \in 
\{0,1,2,3\}$, and so $m_{ij} \in \{2,3,4,6\}$.   
It is known (see for example Proposition 1.3.21 of \cite{Kumar}) 
that the associated {\em Weyl group} is $W(\Gamma,M)$.  
%

\noindent 
{\bf \ProofCorollary}\ \ {\sl (1) Given a generalized Cartan matrix,   
the associated Weyl group is finite if and only if the 
associated Kac--Moody  
algebra is finite-dimensional if and only if 
the associated GCM graph is a Dynkin diagram of finite type. 
(2) Given an E-GCM, the associated Coxeter group is finite if and only 
if the associated E-GCM graph is an E-Coxeter graph.}

{\em Proof.} 
For (2) consider an E-GCM graph $(\Gamma,M)$.  If $W(\Gamma,M)$ is 
finite, then it follows from 
\ErikssonWordProposition.1 above that every game sequence played from some 
given strongly dominant position will converge.  Then $(\Gamma,M)$ is 
admissible and by \SecondMainResult\ must be an E-Coxeter graph.  
Conversely, suppose $(\Gamma,M)$ is an E-Coxeter graph.  Pick a 
strongly dominant position $\lambda$.  Since every game sequence 
played from $\lambda$ converges to 
the same terminal position in the same finite number of steps, then 
by \ErikssonWordProposition.2 we have an upper bound on the lengths 
of reduced expressions in $W(\Gamma,M)$.  Then $W(\Gamma,M)$ is 
finite.  
For (1) observe that 
this same reasoning together with \MarsAttacksTheorem\ shows 
that for a GCM graph $(\Gamma,M)$, the associated Weyl group 
$W(\Gamma,M)$ is finite if and only if $(\Gamma,M)$ is a Dynkin 
diagram of finite type.  By the root space decomposition of 
$\mathfrak{g}(\Gamma,M)$ (\cite{Kumar} Theorem 
1.2.1) and Proposition 1.4.2 of \cite{Kumar}, we have that 
$\mathfrak{g}(\Gamma,M)$ is finite-dimensional if and only if 
$W(\Gamma,M)$ is finite.\hfill\QED


It is well known that the Kac--Moody  
algebras associated to the Dynkin 
diagrams of \DynkinDiagramFigure\ are the complex 
finite-dimensional simple Lie algebras (see for example \cite{Hum} \S 
18). 
It is also well known that Lie algebras corresponding to distinct 
Dynkin diagrams of \DynkinDiagramFigure\ are non-isomorphic.    
For the associated Weyl groups, the only 
redundancy is that the 
groups corresponding to the $\myB_{n}$ and $\myC_{n}$ graphs for $n \geq 3$  
are the same. 
The irreducible Coxeter groups associated to two connected E-Coxeter 
graphs are isomorphic if and only if the graphs are of the same 
type $\mathcal{X}_{n} \in \{\mathcal{A}_{n}, \mathcal{B}_{n}, 
\mathcal{D}_{n}, \mathcal{E}_{6}, \mathcal{E}_{7}, \mathcal{E}_{8}, 
\mathcal{F}_{4}, \mathcal{H}_{3}, \mathcal{H}_{4}, 
\mathcal{I}_{2}(m)\}$.

\vspace{1ex} 

\noindent
{\Large \bf \RemarksNum.\ \ Comments on connections 
with other work}

\vspace{1ex} 
Say an E-GCM graph is {\em strongly admissible} if every nonzero 
dominant position has a convergent game sequence.  
In \cite{ErikssonThesis}, it is shown that:  
{\sl A connected E-GCM graph 
is strongly admissible if and only if it is a connected E-Coxeter 
graph.} This statement essentially 
combines Theorems 6.5 and 6.7 of \cite{ErikssonThesis}.  
Wildberger re-derives this result in \cite{WildbergerPreprint}.    
See \cite{ErikssonLinear} for an ``A-D-E'' version.  

For a Cartan matrix $M$, the 
$M$-structure property of \S \RankedPosetsNum\ is a necessary 
condition for an edge-colored ranked poset to carry certain 
information about some finite-dimensional 
representation of the corresponding semisimple Lie 
algebra.  
Indeed, identifying such combinatorial properties is part of a program 
described in \cite{DonSupp} for 
obtaining combinatorial models for Lie algebra representations.  For 
example, in \cite{ADLMPPW}, four families of finite 
edge-colored distributive lattices are introduced, one family for each 
of the four rank two semisimple Lie algebras.  
Each lattice possesses the $M$-structure property for a Cartan 
matrix $M$ corresponding to the appropriate rank two semisimple Lie 
algebra.  
The ``weight-generating functions'' on these lattices are Weyl 
characters for the irreducible representations of the rank two 
semisimple Lie algebras.  In \cite{DonTwoColor}, the 
posets of join irreducibles (cf.\ \cite{Stanley}) for these 
distributive lattices are shown to be characterized by 
a short list of combinatorial properties.  These are called 
``semistandard posets'' in \cite{ADLMPPW}; the smallest of these posets 
are called ``fundamental posets.''  In \cite{DW}, we will say how 
these fundamental and semistandard posets can be constructed from 
information obtained by playing the numbers game on two-node Dynkin 
diagrams of finite type.  More generally, for an $n$-node Dynkin diagram  
of finite type we will show how to construct other fundamental 
posets from numbers games played from 
certain special fundamental positions, namely the ``adjacency-free'' 
fundamental positions classified in \cite{DonEnumbers}.  
These fundamental posets can be combined to obtain 
semistandard posets whose corresponding distributive lattices produce 
the Weyl characters (as above) for certain irreducible 
representations of the corresponding semisimple Lie algebra.  
It is natural to ask to what extent such distributive lattice models 
for Weyl characters can be characterized combinatorially as in 
\cite{DonTwoColor}. 

\noindent 
{\bf Acknowledgments}\ \ 
We thank  
John Eveland for stimulating discussions during his work on 
an undergraduate research project \cite{Eveland} at Murray 
State University; this helped lead to the questions answered by our 
main results.  
We thank 
Norman Wildberger for sharing his perspective on the numbers game; 
this helped us formulate the question as well as 
our proof of our first main result.  
We thank Bob Proctor for many helpful 
communications during the preparation of 
this paper.

\vspace*{-0.2in}

\renewcommand{\baselinestretch}{1}

\end{document}